\documentstyle {amsppt}
\magnification=1200
\input psfig.sty
\refstyle{C}
\NoBlackBoxes
\TagsOnRight
\NoRunningHeads
\tolerance=3000
\nologo
\voffset=0pt
\hoffset=0pt
\pagewidth{6.5 true in} \pageheight{9 true in}
\def\SBIMSMark#1#2#3{
 \font\SBF=cmss10 at 10 true pt
 \font\SBI=cmssi10 at 10 true pt
 \setbox0=\hbox{\SBF Stony Brook IMS Preprint \##1}
 \setbox2=\hbox to \wd0{\hfil \SBI #2}
 \setbox4=\hbox to \wd0{\hfil \SBI #3}
 \setbox6=\hbox to \wd0{\hss
             \vbox{\hsize=\wd0 \parskip=0pt \baselineskip=10 true pt
                   \copy0 \break%
                   \copy2 \break%
                   \copy4 \break}}
 \dimen0=\ht6   \advance\dimen0 by \vsize \advance\dimen0 by 8 true pt
                \advance\dimen0 by -\pagetotal
 \dimen2=\hsize \advance\dimen2 by .25 true in
%
%
  \openin2=publishd.tex
  \ifeof2\setbox0=\hbox to 0pt{}
  \else 
     \setbox0=\hbox to 3.1 true in{
                \vbox to \ht6{\hsize=3 true in \parskip=0pt  \noindent  
                \input publishd.tex 
                \vfill}}
  \fi
  \closein2
  \ht0=0pt \dp0=0pt
 \ht6=0pt \dp6=0pt
 \setbox8=\vbox to \dimen0{\vfill \hbox to \dimen2{\copy0 \hss \copy6}}
 \ht8=0pt \dp8=0pt \wd8=0pt
 \copy8
 \message{*** Stony Brook IMS Preprint #1, #2 ***}
}

\SBIMSMark{1996/12}{October 1996}{}
\openup 1pt

\bigskip
\bigskip
\bigskip
\bigskip
\topmatter
\title
Universal Models for Lorenz Maps
\endtitle
\author
 Marco Martens \footnote {Institute of Mathematical Sciences, SUNY at 
Stony Brook, Stony Brook, NY 11794-3651.} 
and Welington de Melo \footnote {IMPA, Estrada Dona Castorina 110, 
22460-320 Rio de Janeiro, Brazil. }
\endauthor
\date 
{September 23, 1996}
\enddate
\endtopmatter
\bigskip
\bigskip
\bigskip
\bigskip
\bigskip
\bigskip
\centerline {\bf Abstract.} 

\flushpar
The existence of smooth families of Lorenz maps exhibiting all possible dynamical behavior is established and the structure of the parameter space of these families is described.
 
\flushpar

\newpage

\bigskip
\centerline{\bf 1. Introduction}
\bigskip
   
\flushpar
The aim of this paper is to exhibit some parameterized families of Lorenz 
flows  that are topologically universal in the sense that given any geometric 
Lorenz flow, its dynamics is essentially the same as the dynamics of some element of the family. Thus, these families plays, in the context of Lorenz flows, the same role as the quadratic family in the context of unimodal interval maps. 

\flushpar
Lorenz in [L] showed numerically the existence of some flows in three dimension that have  complicated recurrent behavior. What we now call a 
Lorenz flow has a singularity of saddle type with a one dimensional unstable 
manifold and an infinite set of hyperbolic periodic orbits, whose closure
contains the saddle point. More specifically, the closure of this set of periodic orbits is in general the global attractor of the flow.

\flushpar
To analyze the dynamics of such a flow we take a two 
dimensional transversal section intersecting the local stable manifold in a 
line $l$ and we look at the first return map to $S$. This map is not defined 
in the line $l$ and in fact exhibits a discontinuity at $l$ because orbits 
near $l$ in opposite sides follows different branches of the unstable manifold. 
To describe the dynamics of such a flow, Guckenheimer and Williams added a new hypothesis: 
the existence of a one dimension foliation in $S$ that is invariant by the 
first return map, has $l$ as a leaf and is such that points in the same leaf 
are exponentially contracted under iteration by the first return map. A Lorenz 
flow with this extra structure we call a {\it geometric Lorenz flow}. Because of the 
exponential contraction on the leaves of the foliation, the dynamics of such a 
flow can be described by the action of the first return map on the 
space of leaves of the stable foliation. This space of leaves is an interval 
and the induced map has a unique discontinuity at the point corresponding to $l$. Such an interval map we call a Lorenz map. More precisely,

\proclaim{Definition 1.1} Let $P<0<Q$. A {\it Lorenz} map from $[P,Q]$ to
$[P,Q]$ is a pair $(f_-,f_+)$ where
\item{1)}
$f_-:[P,0]\to [P,Q]$ and $f_+:[0,Q]\to [P,Q]$
are continuous and strictly increasing maps. 
\item{2)} $f(P)=P$ and $f(Q)=Q$.
\item{3)} Given $\rho>0$, we will say that $f$ is a $C^r$ of exponent $\rho$
if we can write 
$$
f_-(x)=\tilde{f}_-(x^\rho) \text{ and } f_+(x)=\tilde{f}_+(x^\rho)
$$
where $\tilde{f}_-$ and $\tilde{f}_+$ are $C^r$ diffeomorphisms defined
on appropriate closed intervals.

\flushpar
Notation: This Lorenz map is denoted by $(P,Q,f_-,f_+)$.
\endproclaim 

\flushpar
 Notice that if $r\geq 1$ then the triple $\{\rho, \tilde f_-,\tilde f_+\}$  is uniquely determined by $f$. If the map is associated to a Lorenz vector field then the exponent $\rho$ is precisely the absolute value of the ratio between the unstable eigenvalue and the weak stable eigenvalue at the saddle point. 

\flushpar
A Lorenz map is {\it non-trivial} iff $f([P,0])\supset [P,0]$ and 
$f([0,Q])\supset [0,Q]$. Otherwise $f$ is trivial, any orbit of such a 
map is asymptotic to a fixed point.  

\flushpar
Guckenheimer and Williams  proved in [GW] that there exists an open set of vector fields 
in three space that have a structure of geometric Lorenz flow with smooth 
associated Lorenz maps. In fact they only considered the situation where 
the exponent is  smaller than one and the map expanding with derivative 
everywhere bigger than $\sqrt 2$. However we can use the same arguments to 
construct open sets of vector fields having Lorenz maps with exponent bigger 
than one. As we will see soon, the  Lorenz maps of exponent bigger 
than one presents a much bigger variety of dynamical behavior due to the interplay of contractions and expansions. Compare this with the unimodal situation: the quadratic family exhibits more types of combinatorics than the expanding tent family.
Before stating our results we need to discuss some combinatorial aspects of 
Lorenz maps. 

\flushpar
A {\it branch of $f^n$} is a maximal closed interval $J$ such that $f^n$ is a diffeomorphism in the interior of $I$.
So an end point of $J$ is either $P$, or $Q$ or a point in the 
backward orbit of $0$. To each branch $J$ of $f^n$ we can associate a word 
$\alpha= (\alpha_0, \alpha_1, \dots, \alpha_{n-1})$ where $\alpha_i\in \{[P,0], 
[0,Q]\}$ and $f^i(J)\subset \alpha_i$. It is clear that given a word of length $n$, there exists at most one branch of $f^n$ associated to it. The combinatorics of all possible words are determined by the kneading invariants of $f$, $K_-(f)$ and 
$K_+(f)$ defined as follows: the first $n$ symbols of $K_-(f)$ are the symbols of the branch of $f^{n+1}$ adjacent to $0$ that is contained in 
$[P,0]$. Similarly for $K_+(f)$. There are many papers describing the combinatorics of Lorenz maps, see, for example, [P] and specially [HS] where 
all possible kneading invariants of Lorenz maps are characterized. 

\flushpar
The intersection of all branches that contains a given point is either a point or 
a closed interval. If such an intersection is an interval, it is called a {\it homterval} of $f$. A {\it critical homterval } is a homterval that has $0$ as an endpoint. So there are at most two critical homtervals. 
The image of a homterval is always contained in another 
homterval and, if the homterval is not critical, it is onto. An orbit of homterval is a sequence $J_0, J^1, \dots $ of homtervals such that $f(J^n) \subset J^{n+1}$. There are three types of orbits of homtervals: 
\parindent=15pt
\item{1)} $J^0$ is a  {\it wandering  interval} if its orbit contains infinitely many intervals; 
\item{2)} it is {\it periodic} of period $n$ if $J^n=J_0$ and  $J^i\neq J^0$ if $0<i<n$; 
\item{3)} it is 
{\it eventually periodic} if it is not periodic but $J^i$ is periodic for some $i$.

\proclaim {Definition 1.2} We say that $f$ is a {\it simple Lorenz map} if
 $f$ has no homterval.
 \endproclaim

\proclaim{Definition 1.3} A {\it maximal  semi conjugacy} from a Lorenz map $f$ to a Lorenz map $\hat f$ is a continuous, monotone, surjective map $h$ 
such that: 
\item{1)}  $\hat f \circ h = h\circ f$ ; 
\item{2)} the inverse image of each point is either a point or a homterval or 
an interval whose points are all asymptotic to periodic points.
\endproclaim

\flushpar
In section 6 we will discuss this notion of semi conjugacy. Roughly speaking,
a maximal semi conjugacy collapses as much as possible without destroying essential parts of the dynamics. If $f$ has at most one critical homterval then a maximal semi conjugacy collapses only homtervals. 

\proclaim{Definition 1.4} The Lorenz maps $f_1$ and $f_2$ are said to be essentially conjugated iff there exists a simple Lorenz map $\hat{f}$ and  nice 
semi conjugacies $h_1$, from $f_1$ to $\hat{f}$, and $h_2$ from $f_2$ to $\hat{f}$. 
\endproclaim

\flushpar
Let $\Cal{L}^r$ be the collection of $C^r$, $r\ge 0$, Lorenz maps. We endow  $\Cal{L}^r$ with a topology that
takes care of the domain (P,Q are close), of the exponents and of the coefficients: the coefficients, after a linear rescaling, to make the domains
equal, are $C^r$ close to each other.

\proclaim{Definition 1.5} Let $\Lambda\subset \Bbb{R}^2$ be closed. A Lorenz
family is a continuous map $F:\Lambda\to \Cal{L}^r$, $r\ge 0$
$$
F_\lambda=\{P_\lambda,Q_\lambda,\phi_\lambda,\psi_\lambda\}.
$$
A Lorenz family $ \lambda \in \Lambda \mapsto F_\lambda$ is {\it full} if given 
any non-trivial Lorenz map $f\in \Cal{L}^2$ there exists a parameter value $\lambda$ 
such that $f$ is essentially conjugated to $F_\lambda$.
\flushpar
A Monotone Lorenz family is a  $C^3$ Lorenz family such that
\parindent=15pt
\item {1)} $F_\lambda$ has negative Schwarzian derivative for all $\lambda\in \Lambda$. 
\item{2)} $\Lambda=[0,1]\times [0,1]$.
\item{3)} $F:(s,t)\to \{-1,1,\phi_s,\psi_t\}$.
\item {4)} if $s_1<s_2$ then $\phi_{s_1}(x)<\phi_{s_2}(x)$ for all 
$x\in [-1,0]$ and if $t_1<t_2$ then $\psi_{t_1}(x)<\psi_{t_2}(x)$ for all 
$x\in [0,1]$
\item{5)} $\phi_0(0)=0$, $\phi_1(0)=1$, $\psi_0(0)=-1$ and $\psi_1(0)=0$. 
\item{6)} $DF_\lambda(\pm 1)>1$ for $\lambda\in \Lambda$.
\endproclaim

\proclaim { Theorem 1.6} A monotone Lorenz family is a full family. 
\endproclaim 

\flushpar
 Next we discuss the notion of  renormalization for Lorenz maps. This will allow
us to refine the above Theorem.

\proclaim {Definition 1.7} A Lorenz map $f$ is called 
renormalizable if there exist $P<p<0<q<Q$ such that the first return map
to $(p,q)$ is a Lorenz map, say $(p,q,f^a,f^b)$. This induced 
Lorenz map is called a renormalization of $f$. The interval $[p,q]$
is called a domain of renormalization, $(a,b)$ are the periods of renormalization and the type of the renormalization is the pair $(\alpha, \beta)$, where $\alpha $ (resp. $\beta$) is the word 
associated to the branch of $f^a$  (resp. $f^b$) that contains $[p,0]$ (resp. $[0,q]$). 
\endproclaim

\flushpar
Let $\Cal{D}_{\alpha,\beta}\subset \Cal{L}^r$ be the subset of Lorenz maps
which have a renormalization of type $(\alpha,\beta)$. These sets are called 
{\it Domains of Renormalization}. In section 3 it will
be shown that 
\item{1)} $\Cal{D}_{\alpha,\beta}$ is closed and connected,
\item{2)} The collection of sets $\Cal{D}_{\alpha,\beta}$ is nested. If
$\Cal{D}_{\alpha,\beta}\cap\Cal{D}_{\hat{\alpha},\hat{\beta}} \ne \emptyset$ 
then 
$$
\Cal{D}_{\alpha,\beta}\subset \Cal{D}_{\hat{\alpha},\hat{\beta}} \text{ or }
\Cal{D}_{\hat{\alpha},\hat{\beta}} \subset \Cal{D}_{\alpha,\beta}.
$$
Let $\Cal{D}\subset \Cal{L}^r$ be the set of renormalizable maps. This set is the union of the nested collection consisting of the sets $\Cal{D}_{\alpha,\beta}$.

\flushpar
Fix a Lorenz family $F:\Lambda\to \Cal{L}^r$. This family will intersect the domains $\Cal{D}_{\alpha,\beta}$, giving rise to the following 

\proclaim{Definition 1.8} The archipelago of type $(\alpha,\beta)$ is the set of are all parameter values $A_{\alpha,\beta}$ for which the corresponding Lorenz map is in $\Cal{D}_{\alpha,\beta}$: $A_{\alpha,\beta}=F^{-1}(\Cal {D}_{\alpha,\beta})$.

\flushpar
An island in the archipelago $A_{\alpha,\beta}$ is a connected component of
the interior of $A_{\alpha,\beta}$. 

\flushpar
To express the type of renormalization we are considering, we will speak about
$(\alpha,\beta)$-archipelagoes and there $(\alpha,\beta)$-islands.
\endproclaim

\flushpar
The archipelagoes inherit properties from the sets $D_{\alpha,\beta}$,
they are closed and nested. This implies that the 
closure of an island defines a Lorenz family, namely the family of the 
corresponding $(\alpha,\beta)-$renormalizations. We call an island a 
{\it full island} if the induced family is a full Lorenz family.

\proclaim{Theorem 1.9} Every archipelago of a monotone Lorenz family contains
a full island. 
\endproclaim   

\flushpar
Observe that Theorem 1.6 is a special case of Theorem 1.9.

\proclaim{Definition 1.10} A Lorenz map is called hyperbolic iff both critical
orbits tend to hyperbolic periodic attractors and the complement of the basin of
these periodic attractors is a hyperbolic set.
\endproclaim

\proclaim{Proposition 1.11} The hyperbolic elements in a monotone Lorenz family form
an open and dense set in parameter space.
\endproclaim

\bigskip
\centerline{\bf Conjectures and Remarks}
\bigskip

\flushpar
We finish this section with some conjectures, problems  and remarks. 

\proclaim {Conjecture 1.12} If $f$ is a $C^2$ Lorenz map with exponent $\rho >1$ then the number of periodic homtervals of $f$ is finite. 
\endproclaim 

\proclaim{Definition 1.13} We say that a Lorenz map $f$ has a {\it Cherry attractor} 
if there exists a renormalization of $f$, $(p,q, f^a,f^b)$ with the following 
properties:
\item 1. The interval $[f^b (0+), f^a (0-)]= [p', q']$ is invariant under the renormalized map  whose restriction $g$ to this interval is one to one but not onto.  
\item 2.  $g$ has no periodic point. 
\flushpar 
The $f-$invariant set $\cup _{i=0}^a f^i([p',0] \cup_{j=0}^qf^j([0,q']$ is 
called a Cherry attractor for $f$.  
\endproclaim
 
\flushpar
If $J$ is the interval $[p',q'] \setminus g([p',q']$ then the inverse of 
$g$ can be extended continuously to the whole interval and gives a map that 
is constant on $J$, strictly monotone otherwise and maps $p',q'$ into $0$. 
Hence $g^{-1}$ can be thought as a circle map with a flat top without periodic 
point. This is called a Cherry map and appears as a first return map of 
a recurrent flow on the torus (see [MMMS]). It follows that $J$ is a wandering 
homterval of $f$.

\proclaim {Conjecture 1.14} Let $f$ be a $C^2$ Lorenz map with exponent $\rho>1$. 
If $f$ has a wandering homterval then $f$ has a Cherry attractor. 
\endproclaim 

\flushpar
In [MMS] it is proved that smooth interval maps have only finitely many periodic
homtervals. Exactly the same proof can be applied to Lorenz maps which do not have wandering intervals. It follows that Conjecture 1.12 is consequence of Conjecture 1.14.

\proclaim { Conjecture 1.15} If $f$ is a $C^2$ Lorenz map with exponent $\rho<1$ then:
\item {1)}   $f$ has at most a finite number of renormalizations.
\item {2)} If $f$ is not renormalizable then either there exists a maximal
semi conjugacy from $f$ to a piece wise affine Lorenz map with constant derivative ( $\beta$-transformation) or 
the restriction of $f$ to the interval $[f(0+),f(0-)]$ is 1-1. 
\endproclaim

\flushpar
Let $\Cal D$ denotes the set of Lorenz maps that are renormalizable. We prove 
in section 3 that each connected component of $\Cal D$ is equal to some 
$\Cal {D}_{\alpha,\beta}$. From this we can define the renormalization operator 
$\Cal {R} \colon \Cal D \to \Cal {L}^r$ as follows. Let $f\in \Cal{D}_{\alpha,\beta}$ and define $\Cal{R}(f) = A\circ \hat f\circ A^{-1}$ where $A(x)= \frac xq$ and $\hat{f}=(p,q,f^a,f^b)$ the $\alpha,\beta-$renormalization of $f$. Hence $\Cal {R}(f)$ is a Lorenz map with the 
positive fixed point equal to $1$. If we restrict $\Cal R$ to the space of 
Lorenz maps with the same normalization we get an operator. 

\proclaim{Conjecture 1.16} Let $\Cal {D}_i, i\in \bold N$ be the connected components of the domain of the renormalization operator. 
\item{1)} Given a finite sequence $i_0, \dots, i_{n-1}$ of integers not necessarily distinct, there exists a unique 
normalized Lorenz map $g$ such that $\Cal R^n(g)=g$ and $\Cal R^k(g) \in 
\Cal {D}_{i_j}$ whenever $k=i_j \mod n$, compare [ACT].

\item{2)} If $f$ is a normalized Lorenz map such that $\Cal R^k(f)\in \Cal D_{i_j}$ whenever $k= i_j \mod n$ 
then $\Cal R^{kn} (f)$ converges exponentially 
fast to
$g$ as $k\to \infty$, compare [ACT] .
 
\item{3)} Let  $f$ be a Lorenz map whose forward orbit under the renormalization operator meets only a finite number of connected components of the domain. 
Then there exists a compact set $\Cal K\subset \Cal L$ such that $\Cal R^n(f)$ belongs to $\Cal K$ for all $n$. In particular, the length of the left component of the domain of $\Cal {R}^n(f)$ is bounded from above and from below. 

\item {4)} Let $\Cal B \subset \Cal D$ be the union of a finite number of 
connected components of $\Cal D$. Then there exists a compact set $\Cal K_\Cal B$ of the space of Lorenz maps such that for any mapping $f\in \Cal B$ whose 
forward orbit remains in $\Cal B$, there exists $n_0$ such that for $n\geq n_0$,

$\Cal R^n(f) \in \Cal K_\Cal B$. 
\endproclaim

\proclaim {Conjecture 1.17}  
There exist a monotone Lorenz family such that each 
archipelago is an island.
\endproclaim 

\proclaim {Conjecture 1.18} For generic Lorenz families the following holds: 
there exists a bound, depending on the family, for the number of island 
in each archipelago. 
\endproclaim 

\flushpar
In [MP] it has been proved that in a $C^1$ generic one-parameter family of $C^2$  circle
difeomorphisms the rotation number is a piecewise monotone function. Conjecture 
1.18 is the corresponding statement for Lorenz families.  

\proclaim { Conjecture 1.19} For a monotone Lorenz family the diameter of the 
islands goes to zero as the period of renormalization goes to infinity.
\endproclaim 

\proclaim { Conjecture 1.20} For a monotone Lorenz family, the set of parameter 
values that belong to infinitely many archipelagoes has Lebesgue measure zero. 
\endproclaim

\flushpar
We say that an  island is of generation 0 if it is not contained in another 
island. By induction we say that an island is of generation $n$ if it is 
not of generation $n-1$ and any island that contains it is of generation $\leq n-1$. From theorem 1.11 it follows that any full island of generation $n$ contains 
infinitely many island of generation $n+1$. In particular there are uncountably many parameter values that corresponds to infinitely renormalizable islands. This is in sharp contrast with Conjecture 1.15. 

\bigskip

\flushpar
From [R] it follows that for generic two parameter families of Lorenz maps, 
the set of parameter values corresponding to maps that have positive Lyapunov 
exponents has positive Lebesgue measure. 
 
\proclaim {Conjecture 1.21} For generic two parameter families of Lorenz maps 
the set of parameter values corresponding to maps that are not hyperbolic and 
that do not have positive Lyapunov exponents has zero Lebesgue measure. 
\endproclaim

\demo{Acknowledgements} Part of this work was written while the first author was visiting IMPA and the second author was visiting IHES and the graduate center
of CUNY. The authors would like to thank these institutions for their kind
hospitality.
\enddemo

\bigskip
\centerline{\bf Notation}
\bigskip

\flushpar
Let $I\subset [-1,1]$ be an interval with boundary points $a$ and $b$, say $a<b$.
Then $\partial I=\{a,b\}$, $\partial_-I=\{a\}$ and $\partial_+I=\{b\}$. If
$J\subset I\subset [-1,1]$ are two intervals and $\partial_-J=\partial_-I$, we
say $J\subset_l I$. If $\partial_+J=\partial_+ I$ then we say $J\subset_r I$.
In the case that $\partial J\cap \partial I=\emptyset$ we write $J\subset_{\text{int}} I$.

\flushpar
The discontinuity of Lorenz maps in $0$ causes that such maps has two
{\it critical orbits}:
$$
f^n(0_-)=\lim_{x\uparrow 0} f^n(x) \text{ and } f^n(0_+)=\lim_{x\downarrow 0} f^n(x), 
$$
with $n\ge 0$.

\bigskip
\centerline{\bf 2. Combinatorial Properties of Lorenz-maps }
\bigskip

\flushpar
We will start by defining kneading sequences for Lorenz maps, similarly as was 
done in [MT] for continuous piecewise monotone interval maps. 
Clearly, to describe Lorenz maps we will need 
two kneading sequences, one for $0_-$ and one for $0_+$. Fix a Lorenz map
$f:[-1,1]\to [-1,1]$.

\bigskip

\flushpar
Let $\Bbb{B}_n(f)$ be the collection of {\it branches of $f^n$}. That is, the collection of maximal intervals on which $f^n$ is monotone. If $I\in \Bbb{B}_n(f)$ then the word $\omega(I)\in \{L,R\}^n$ is such that
$$
\aligned
\omega_i(I)&=L \text{ if } f^i(I)\subset [-1,0)\\
\omega_i(I)&=R \text{ if } f^i(I)\subset (0,1]\\
\endaligned
$$ 
for $0\le i<n$. The {\it kneading sequences} are defined as
$$
K^-_n(f)=\omega(I_-(n)) \text{ and } K^+_n(f)=\omega(I_+(n))
$$
where $I_-(n), I_+(n)\in \Bbb{B}_n(f)$, $n\ge 0$, with $\partial_+ I_-(n)=0$
and $\partial_- I_+(n)=0$. Let $K^-(f)$ and $K^+(f)$ be the limits of 
respectively $K^-_n(f)$ and $K^+_n(f)$, $n\to\infty$.

\bigskip

\flushpar
For every branch $I\in \Bbb{B}_n(f)$ there exist unique {\it cutting
times} $l_n(I,f)$ and $r_n(I,f)$ such that 
$$
0\in \partial_+f^{l_n(I,f)}(I) \text{ and } 0\in \partial_-f^{r_n(I,f)}(I).
$$
When it is clear which map is under consideration, we will suppress
 the symbol $f$ in the above notation. 

\proclaim{Lemma 2.1} Let $f$ and $g$ be Lorenz maps such that $K_n^\pm(f)=K_n^\pm(g)$ for some $n\ge 1$ and assume that 
$$
f^i(0_\pm)=0 \Leftrightarrow g^i(0_\pm)=0,
$$
where $i\le n$.

\flushpar
There exists an orientation
preserving homeomorphism $h:[-1,1]\to [-1,1]$ such that
\parindent=15pt
\item{1)} $h$ preserves branches and their type: for every $k\le n$, $I\in \Bbb{B}_k(f)$ 
$$
h(I)\in \Bbb{B}_k(g) \text{ and } \omega(I,f)=\omega(h(I),g).
$$
\item{2)} if $I\in \Bbb{B}_k(f)$ and $f(I)=I'\in \Bbb{B}_{k-1}(f)$ then
$$
g(h(I))=h(I')\in \Bbb{B}_{k-1}(g).
$$
\item{3)} if $I\in \Bbb{B}_k(f)$ and $f(I)\subset_{l/r} I'\in \Bbb{B}_{k-1}(f)$ then
$$
g(h(I))\subset_{l/r} h(I')\in \Bbb{B}_{k-1}(g).
$$
\item{4)} for $k\le n$, $I\in \Bbb{B}_k(f)$ 
$$
l_k(I,f)=l_k(h(I),g) \text{ and } r_k(I,f)=r_k(h(I),g). 
$$
\endproclaim

\demo{Proof} First we will show that 2), 3) and 4) follow from 1). 

\demo{ $1) \Rightarrow 2), 3)$} Let $I=(x,y)\in \Bbb{B}_k(f)$ and 
assume $f(I)\subset I'=(x',y')\in\Bbb{B}_{k-1}(f)$. It suffices to prove

\proclaim{Claim} 
$$
\aligned
f(x)=x' &\Leftrightarrow g(h(x))=h(x') \\ 
f(y)=y' &\Leftrightarrow g(h(y))=h(y') 
\endaligned
$$
\endproclaim

\flushpar
By symmetry we only have to consider the left boundary.
First observe that $f(x)=x'$ whenever $x\ne 0$. So if $x\ne 0$ then 
$h(x)\ne 0$ and $g(h(x))=h(x')$. The Claim is proved in this case.

\flushpar
Assume $x=0$ and $f(x)=x'$. Because $x'$ is a boundary point of $I'\in \Bbb{B}_{k-1}(f)$ there exists $j\le k-1$ such that $f^j(0_+)=0$. By assumption, the same holds for $g$:  $g^j(0_+)=0$. Now assume by contradiction that
$g(h(x))\ne h(x')$. Then $g(0_+)=g(h(x))$ belongs to the interior of $g(h(I'))$ and $g^{j-1}(h(I'))$ contains $0$ in its interior. Therefore, 
$g^j|h(I')$ is not monotone. This is a contradiction because $h(I')\in \Bbb{B}_{k-1}(g)$ and $j\le k-1$. We proved that $g(h(x))=h(x')$. 
\hfill\hfill\qed $\,\,$ ( $1) \Rightarrow 2), 3)$)
\enddemo

\demo{ $1) \Rightarrow 4)$} Observe that for each $I\in \Bbb{B}_k(f)$ and 
$i\le k$ there exists $I'\in \Bbb{B}_{k-i}(f)$ such that one of the following three possibilities holds:
\item{a)} $f^i(I)=I'$,
\item{b)} $f^i(I)\subset_l I'$ or $f^i(I)\subset_r I'$,
\item{c)} $f^i(I)\subset_{int} I'$.

\flushpar
Moreover observe that there are unique numbers $0\le i_1< i_2<k$ such that
for $0\le i <i_1$ $f^i(I)$ is in case a). For $i_1\le i<i_2$ $f^i(I)$ is in case b), say $f^i(I)\subset_l I'$. Finally for $i\ge i_2$ $f^i(I)$ is in case c). Clearly this numbers are exactly the cutting times 
$$
\aligned
i_1&=r_k(I,f) \text{ and } i_2=l_k(I,f) \text{ or }\\
i_1&=l_k(I,f) \text{ and } i_2=r_k(I,f). 
\endaligned
$$
Using properties 2) and 3) of the Lemma it follows that the intervals 
$g^i(h(I))$ jump exactly at the same times $i_1$ and $i_2$ from case to case.
Hence 
$$
\aligned
r_k(h(I),g)&=i_1=r_k(I,f) \text{ and } l_k(h(I),g)=i_2=l_k(I,f) \text{ or }\\
r_k(h(I),g)&=i_2=r_k(I,f) \text{ and } l_k(h(I),g)=i_1=l_k(I,f).
\endaligned
$$
\hfill\hfill\qed $\,\,$ ( $1) \Rightarrow 4)$)
\enddemo

\flushpar
The proof of property 1) will be by induction in $n\ge 1$. It clearly holds for $n=1$. Assume 1) holds for some $n\ge 1$. In particular we may assume that also
2), 3) and 4) hold for this $n\ge 1$. 

\flushpar

Observe that the collection of branches $\Bbb{B}_k(f)$ with $k=0,1,2,\dots ,n$ 
define a refining sequence of partitions of $[-1,1]$. The homeomorphism $h$ 
maps this sequence of partitions to the partitions formed by $\Bbb{B}_k(g)$ with
$k\le n$. To prove 1) for $n+1$ we have to construct a homeomorphism $H$ which
also preserves the above partitions and moreover maps $\Bbb{B}_{n+1}(f)$ to
$\Bbb{B}_{n+1}(g)$. In particular $H$ is obtained by redefining $h$ in the 
interior of the branches $\Bbb{B}_n(f)$.

\flushpar
Consider a branch $I\in \Bbb{B}_n(f)$ with $h(I)=I'\in \Bbb{B}_n(g)$. The new 
homeomorphism $H$ will also satisfy $H(I)=I'$. The boundary of the interval
$f^n(I)$ consists of the critical values 
$$
\partial f^n(I)=\{ f^{n-l_n(I,f)}(0_+), f^{n-r_n(I,f)}(0_-)\}.
$$
Property 4) of the Lemma states that also
$$
\partial g^n(I')=\{ g^{n-l_n(I,f)}(0_+), g^{n-r_n(I,f)}(0_-)\},
$$
the cutting times of $I$ and $I'$ are the same.

\flushpar
Now we will use that the kneading sequences of $f$ and $g$ are the same up to 
$n+1$. Let $(\theta_i)_{i=1,\dots, n+1}= K^\pm_n(f)=K^\pm_n(g)$.
Observe that
$$
I\notin \Bbb{B}_{n+1}(f) \Longleftrightarrow
\theta_{n-r_n(I,f)}=R \text{ and } \theta_{n-l_n(I,f)}=L.
$$
And 
$$
I'\notin \Bbb{B}_{n+1}(g) \Longleftrightarrow
\theta_{n-r_n(I,g)}=R \text{ and } \theta_{n-l_n(I,g)}=L.
$$
In particular we have 
$$
I\in\Bbb{B}_{n+1}(f)\Longleftrightarrow I'\in\Bbb{B}_{n+1}(g).
$$ 
If $I\in \Bbb{B}_{n+1}(f)$ we do not have to change $h$: $H|I=h$. If $I\notin \Bbb{B}_{n+1}(f)$ then both $I$ and $I'$ will have two branches of respectively
$\Bbb{B}_{n+1}(f)$ and $\Bbb{B}_{n+1}(g)$. Define $H|I$ such that these two
sub-branches are matched.

\flushpar
Once this construction has been done for all branches $I\in \Bbb{B}_n(f)$ we will obtain a homeomorphism $H$ which preserves the branches in $\Bbb{B}_{n+1}(f)$ and their types.

\hfill\hfill\qed $\,\,$ (Lemma 2.1)
\enddemo

\demo{Example} Let $\phi_\lambda:[-1,1]\to [-1,1]$, $\lambda\in [0,1]$, 
be a one parameter
family of unimodal maps. Say, $\phi_\lambda(\pm 1)=-1$ for $\lambda\in [0,1]$, and $\phi_0(0)=0$, $\phi_1(0)=1$. Moreover, assume that this unimodal family is monotone:  whenever $\lambda_1<\lambda_2$ we have $\phi_{\lambda_1}(t)< \phi_{\lambda_2}(t)$ for all $t\in (-1,1)$. 

\flushpar
Consider the Lorenz family $f_{x,y}:[-1,1]\to [-1,1]$ with $(x,y)\in [0,1]\times [0,1]$ and  
$$
\aligned
f_{(x,y)}(t)&=\phi_x(t) \text{ if } t<0\\
f_{(x,y)}(t)&=\phi_{1-y}(t) \text{ if } t>0.
\endaligned
$$
Notice that for $(x,y)\in U=\{(t,s)\in [0,1]\times [0,1]| s=1-t\}$ the Lorenz map $f_{(x,y)}$ behaves exactly like the $U${\it nimodal} map $\phi_x$. Namely
$$
f^n_{(x,y)}(t)=\pm\phi^n_x(t),
$$
for $n\ge 0$ and $t\in [-1,1]$.
\enddemo

\flushpar
This example inspires to change the coordinates of the parameter space of monotone Lorenz families. Fix a monotone Lorenz family $F:[0,1]\times [0,1]\to \Cal{L}$. The usual coordinates $x$ and $y$ for $[0,1]\times [0,1]$ are not appropriate to explore the
monotonicity of the family. We will use 
$$
\aligned
U&=\{(x,y)|y=1-x\}\\
M&=\{(x,y)|y=x\}.
\endaligned
$$
as coordinate axis for $[0,1]\times [0,1]$ giving rise to the coordinates $u$ and $m$ on $[0,1]\times [0,1]$, $u=x-y, m= x+y-1$.

\flushpar
For $u\in U$ let $M_u=\{(u,m)|m\in M\}=\{(x,y)| x-y=u\}$ and for 
$(x,y)\in [0,1]\times [0,1]$ let 
$$
\aligned
C^+_{(x,y)}&=\{(t,s)| t\ge x, s\ge y\},\\ 
C^-_{(x,y)}&=\{(t,s)| t\le x, s\le y\},\\
B_{(x,y)}  &=[0,1]\times [0,1] \setminus ( C^+_{(x,y)}\cup C^-_{(x,y)} ).
\endaligned
$$

\bigskip
\centerline{\psfig{figure=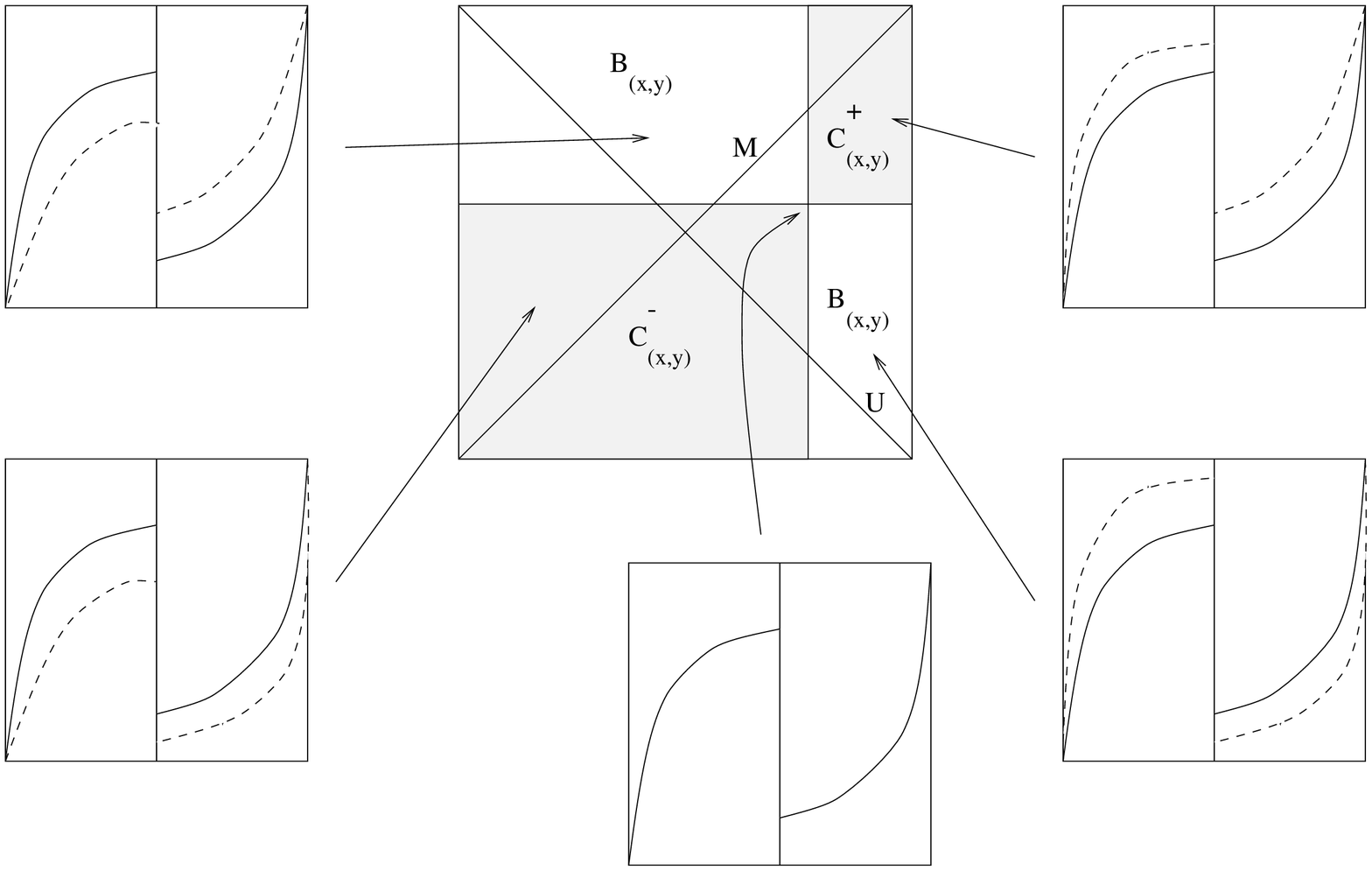,width=12cm}}
\smallskip
\centerline{\bf Figure 1 The $u,m-$coordinates}
\bigskip

\flushpar
These cones will play a crucial role in the study of the parameter space of monotone
Lorenz families. The $C^+-$cones describe deformations in which both branches {\it move up}. The $C^--$cones describe deformations in which both branches {\it move down}. The $B-$cone contains the maps for which one branch moves up and the other down. The situation is illustrated in Figure 1. The deformations into the $B-$cones are the ones which are difficults to understand. We will explore the deformations into the $C^\pm-$cones.

\flushpar
For example, the monotonicity of the family implies immediately the monotonicity in kneading information: let $K^\pm(z)=K^\pm_{F(z)}$ then
$$
K^\pm(z')\ge K^\pm(z) \text{ for } z'\in C^+(z)
$$ 
and
$$
K^\pm(z')\le K^\pm(z) \text{ for } z'\in C^-(z),
$$
where we used the usual lexicographic order on the $L-R$-sequences ( $L<R$).

\demo{Proof of Proposition 1.11} Given a map $f$ in a monotone Lorenz family
we will construct arbitrarily close a map in the family whose critical orbits tend to hyperbolic periodic attractors. Then general arguments ([M])
will show that this perturbation is actually a hyperbolic Lorenz map.

\flushpar
Assume that both critical orbits do not tend to a periodic attractor. We also
may assume that the map $f$ is not in the boundary of parameter space: 
$C^-_f\cap ([0,1]\times [0,1])\ne \emptyset$. The first step will be to find a map in $C^-_f$ close to $f$ such that $0_-$ is periodic.

\bigskip

\flushpar
Take $g_1\in C^-_f$ close to $f$. Assume that there is some minimal $n_0>0$ such that $g_1^{n_0}(0)<0<f^{n_0}(0_-)$. If we move back along the straight line from $g_1$ to $f$ then $g^{n_0}(0_-)$ will increase up to $f^{n_0}(0_-)$. Hence along this straight line there is a map $g$ with $g^{n_0}(0_-)=0$.

\flushpar
The second case is that for all $n>0$ we have $g^n_1(0_-)$ and $f^n(0_-)$ are on the same side of $0$. In particular $g^n_1(0_-)<f^n(0_-)$ for all $n>0$.
Consider the intervals
$$
J_n=[g_1^n(0_-),f^n(0_-)]
$$
with $n\ge 0$. Observe $f(J_n)\subset J_{n+1}$. This interval $J_n$ has to accumulate at $0_+$. If not the interval $J_1$ would be a hometerval of $f$ not accumulating on $0_+$. Then we can define a continuous $C^2$ map, by redefining $f$ on a small interval $(0,a]$,
having $J_1$ as a hometerval: by [MMS] we know that $J_1$ tends to a periodic attractor of $f$. In particular the orbit $f^n(0_-)$ tends to a periodic attractor. This is a contradiction because we assumed both critical orbits not to accumulate at periodic attractors.

\flushpar
We proved that the critical orbit $g_1^n(0_-)$ accumulates at $0_+$. For each
perturbation $g_2\in C^-_{g_1}$ there will be some $n_0>0$ such that 
$g^{n_0}_2(0_-)<0<g^{n_0}_1(0_-)$. Again we can move back a little bit to find a map $g\in C_f^-$ close to $f$ such that $0_-$ is periodic for $g$.
We finished the proof of

\proclaim{Claim} If both critical orbits of $f$ are not accumulating at a periodic attractor then arbitrarily close there is a map $g$ such that $0_-$ is periodic. In particular arbitrarily close there is a map such that $0_-$ is attracted to a hyperbolic periodic attractor but is not itself a periodic attractor.
\endproclaim   

\flushpar
Assume that one critical orbit is attracted to a periodic attractor but is not itself periodic. We may assume that this orbit is the of the critical point  $0_-$.
The boundary point $-1$ and $1$ are expanding fixed points. Because the map has
negative Schwarzian derivative the periodic attractor is of one of the following types.
\item{1)} There exist an interval $(0,p]$ and $n>0$ such that $f^n((0,p))\subset (0,p)$ and $f^n(p)=p$. By taking $g\in C^-_f$ close to $f$ we may assume that the periodic attractor is hyperbolic and still attracts also $0_-$. We got a situation in which both critical orbits are attracted to a hyperbolic periodic attractor. The Proposition 1.11 is proved.

\item{2)} There exist an interval $[p,0)$ and $n>0$ such that $f^n([p,0))\subset [p,0)$. Again by taking $g\in C^+_f$ we may assume that the periodic orbit of $p$ is a hyperbolic attractor. This situation will persist in a small neighborhood of $g$.

\flushpar
Left is to deform $g$ such that also $0_+$ is attracted towards a hyperbolic attractor. Assume that $0_+$ is not attracted towards a hyperbolic attractor.
It could be attracted towards a neutral periodic attractor. Take 
$g_1\in C^+_{g}$ close to $g$. Consider the two orbits $g^n(0_+)$ and $g_1^n(0_+)$, $n\ge 1$. If at some moment these points are separated then we can move back a little bit towards $g$ and make $0_+$ periodic and by moving back a little bit more we can make it to be attracted towards the periodic attractor which already attracted $0_-$.

\flushpar
If these points are never separated we show, as before, that they accumulate at $0_+$ or is attracted towards a periodic attractor. As before we can make $0_+$ to be attracted towards a hyperbolic attractor. 
\hfill\hfill\qed $\,\,$ (Proposition 1.11)
\enddemo

\flushpar
The proof of Proposition 1.11 also shows that the hyperbolic Lorenz map are dense in the whole space $\Cal{L}^r$ of Lorenz maps (this can also be obtained
from the arguments of [R]). Notice that this space contains a closed subspace $\Cal{L}_s$ of symmetric Lorenz maps. These are 
Lorenz maps $f$ with branches $f_-$ and $f_+$ such that $f_-(x)=-f_+(-x)$. To each
symmetric Lorenz map we can associate a unimodal map $\phi_f$, namely
$$
\aligned
\phi_f(x)&=f_-(x) \text{ for } x<0\\
\phi_f(x)&=-f_+(x) \text{ for } x\ge 0.
\endaligned
$$
Observe that $\phi^n(x)=\pm f^n(x)$, for $n\ge 0$. The Lorenz map $f$ can be seen as an orientation preserving "lift" of the unimodal map $\phi_f$. It follows that the density of the hyperbolic Lorenz maps in $\Cal{L}_s$ would imply the density of hyperbolic maps in the space of unimodal maps. This is an
important open question.

\flushpar
One should compare this with the situation of flows on non-orientable surfaces.
These flows are covered by flows on the orientable cover of the manifold. Proposition 1.11 should be compared with the Peixoto Theorem for flows on orientable surfaces.

\vfill\eject
\bigskip
\centerline{\bf 3. Domains of Renormalization}
\bigskip

\flushpar
In this section we will study some topological properties of the sets $\Cal{D}_{\alpha,\beta}$, the Lorenz maps which have a $(\alpha,\beta)-$renormalization.

\proclaim{Lemma 3.1} Let $p<0<q$ be two periodic points of the 
Lorenz map $f$. Say with period $a$ and $b$ respectively.
Assume
\parindent=15pt
\item{1)} $(p,q)\cap(orb(p)\cup orb(q))=\emptyset$.
\item{2)} $f^a|[p,0]$ and $f^b|[0,q]$ are monotone.
\item{3)} $f^a([p,0])\supset [p,0]$ and $f^b([0,q])\supset [0,q]$. 

\flushpar
Let $(l,0)$ (resp. $(0,r)$) be the maximal interval on which $f^a$ 
(resp $f^b$) is monotone.
\flushpar
Then
$$
f^a(l)\le l \text{ and } f^b(r)\ge r.
$$
Furthermore if $f^a(l)=l$ (resp. $f^b(r)=r$) then $f^b(0_+)=0$ 
(resp. $f^a(0_-)=0$). 

\flushpar
In particular, if the Lorenz map has negative Schwarzian derivative then
 the periodic points $p$ and $q$ are hyperbolic repellors.
\endproclaim

\demo{proof} Assume that $f^a(l)\ge l$. Let $L=(l,p)$. The assumption implies that $f^a(L)\subset L$. Consequently we get that $f^n|_L$ is monotone for all 
$n\ge 0$, $L$ is a homterval.

\flushpar
The maximality of the interval $(l,0)$ gives some $i< a$ such that $0\in \partial f^i(L)$. From property 1, it follows that  $f^i(L)\supset (0,q)$. In particular $(0,q)$ is a homterval. This is only possible if $f^b(0_+)=0$, the point $0_+$ is periodic. Hence $f^a(L)=L$ and therefore,  $f^a(l)=l$.

\flushpar
Let $L_-=L=(l,p)$ and $L_+=(p,0)$. We have $f^a(L_\pm)\supset L_\pm$. The Maximal Principle for maps with negative Schwarzian derivative, see [MS],
 implies that $p$ is a repelor.

\hfill\hfill\qed $\,\,$ (Lemma 3.1)
\enddemo

\proclaim{Lemma 3.2} If $f\in \Cal{D}_{\alpha,\beta}\cap \Cal{D}_{\hat{\alpha},\hat{\beta}}$
then $\hat{\alpha}$ and $\hat{\beta}$ are formed by concatenating the words 
$\alpha, \beta$. In particular
$$
\aligned
\hat{\alpha}=\alpha\beta\dots\\
\hat{\beta}=\beta\alpha\dots
\endaligned
$$
and $|\hat{\alpha}|, |\hat{\beta}|\ge |\alpha|+|\beta|$.

\flushpar
Or vice versa, $\alpha$ and $\beta$ can be expressed in terms of $\hat{\alpha}$ and $\hat{\beta}$.
\endproclaim

\demo{Proof} Let $f\in \Cal{D}_{\alpha,\beta}\cap \Cal{D}_{\hat{\alpha},\hat{\beta}}$ with renormalizations
$(p,q,f^a,f^b)$ and $(\hat{p},\hat{q},f^{\hat{a}},f^{\hat{b}})$ of type
respectively $(\alpha, \beta)$ and $(\hat{\alpha}, \hat{\beta})$. We will first show that 
$$
(\hat{p},\hat{q})\subset (p,q) \text{ ( or } (p,q)\subset (\hat{p},\hat{q}) \text{ ). }
$$
Without loss of generality we may assume that $(p,q)$ is not contained in $(\hat{p}, \hat{q})$, say $\hat{p}\in (p,0)$. We are going to show that $\hat{q}\in (0,q)$, the intervals are nested.

\flushpar
Assume $\hat{q}>q$. Then 
$$
f^{\hat{b}}(q)\in f^{\hat{b}}((0,\hat{q}))\subset (\hat{p},\hat{q})\subset (p,\hat{q}).
$$
However, the orbit of $q$ never enters $(p,q)$. So $f^{\hat{b}}(q)\in [q,\hat{q})$. Because the orbit of $q$ is finite it is
impossible that $f^{\hat{b}}(q)\in (q,\hat{q})$: the map
$f^{\hat{b}}:[q,\hat{q}]\to [q, \hat{q}]$ is monotone and would have a periodic attractor.  So $f^{\hat{b}}(q)=q$ and $\hat{b}=m\cdot b$. The interval
$[q,\hat{q}]$ is periodic which is impossible because the orbits in the boundary
are expanding, see Lemma 3.1. We proved $\hat{q}\le q$.

\bigskip

\flushpar
The aim is to prove that $\hat{q}<q$. Assume $\hat{q}=q$. Then both renormalizations will have $f^b|_{(0,q)}$ as right branch. In particular
$f^b(0_+)> \hat{p}$. 

\flushpar
Let $(x,0)\subset (p,0)$ be the maximal monotone interval on which $f^{\hat{a}}$
is monotone. The map $f$ admits a renormalization of type $(\hat{\alpha},\hat{\beta})=(\hat{\alpha},\beta)$ we have 
\item{1)} $\hat{p}\in (x,0)$,
\item{2)} $f^{\hat{a}}((x,0))\supset (x,0)$ (by Lemma 3.1). 

\flushpar
Observe that $f^a((x,0))\subset (0,q)$. The next images of this interval will 
always stay on the right side of $f^b(0_+)$. We have 
$$
f^b(0_+)>\hat{p}>x,
$$
the images of $(x,0)$ will never be able to cover $(x,0)$ completely, contradiction.
We proved $\hat{q}<q$. In particular that the $(\hat{\alpha},\hat{\beta})$
renormalization is a first return map of the $(\alpha,\beta)$ renormalization.
\hfill\hfill\qed $\,\,$ (Lemma 3.2)
\enddemo

\flushpar
In the proof of the next Proposition we will use the usual lexicographic
order on the $0-1-$words. The length of a word $\omega$ will be denoted by 
$|\omega|$. Moreover, if $j<\min\{|\omega_1|,|\omega_2|\}$ then
$$
\omega_1<_j \omega_2
$$
means that $\omega_1< \omega_2$ and the words differ in the first $j$ symbols.

\proclaim{Proposition 3.3} The domains of renormalization have the following properties
\item{1)} $\Cal{D}_{\alpha,\beta}$ is closed and connected,
\item{2)} The collection of sets $\Cal{D}_{\alpha,\beta}$ is nested. If
$\Cal{D}_{\alpha,\beta}\cap\Cal{D}_{\hat{\alpha},\hat{\beta}} \ne \emptyset$ 
then 
$$
\Cal{D}_{\alpha,\beta}\subset \Cal{D}_{\hat{\alpha},\hat{\beta}} \text{ or }
\Cal{D}_{\hat{\alpha},\hat{\beta}} \subset \Cal{D}_{\alpha,\beta}.
$$
\endproclaim

\demo{proof} The sets $\Cal{D}_{\alpha,\beta}$ are closed because of Lemma 3.1. Let $f,g\in \Cal{D}_{\hat{\alpha},\hat{\beta}}$ and assume that $f\in \Cal{D}_{\alpha,\beta}$. To prove the Proposition we have to show that also $g\in \Cal{D}_{\alpha,\beta}$.. Observe that,
by Lemma 3.2, 
the kneading sequences of $f$ and $g$ are equal up to at least $|\alpha|+|\beta|=a+b$.

\bigskip

\flushpar
Let $(p,q,f^a,f^b)$ be the $(\alpha,\beta)$ renormalization of $f$ and 
$(0,y)$ the maximal interval on which $f^b$ is monotone. Lemma 3.1
states that $f^b((0,y))\supset (0,y)$. In particular $f^b(y)>y$

\flushpar
There exists a unique $j=r_b((0,y),f)<b$, which by Lemma 2.1 depends only on $\alpha$ and 
$\beta$, such that $0=\partial_+ f^{j}((0,y))$. The kneading sequence of 
$f^b(y)$ equals $\sigma^{b-j}(\hat{\alpha})$ and satisfies
$$
\sigma^{b-j}(\hat{\alpha})>_j \beta,
$$
which is the combinatorial formulation of $f^b(y)> y$. Observe that
$$
|\sigma^{b-j}(\hat{\alpha})|=\hat{a}-(b-j)\ge a+b-(b-j)>j,
$$
the above relation between $\sigma^{b-j}(\hat{\alpha})$ and $\beta$ is well defined. 

\flushpar
The map $g$ has a 
$(\hat{\alpha},\hat{\beta})$ renormalization. Lemma 3.2 states that
$\hat{\beta}=\beta\alpha\cdots$. The map $g$ has a branch $(0,y')$ of $g^b$
of type $\beta$. We are going to prove that 
$$
g^b((0,y'))\supset (0,y').
$$
First, because $\hat{\beta}=\beta\alpha$ we have $g^b(0_+)<0$. The branch $(0,y')$ is going to be cut also at moment $j=r_b((0,y'),g)=r_b((0,y),f)$. So the kneading sequence of $g^b(y')$ equals
$$
\sigma^{b-j}(\hat{\alpha})>_j \beta.
$$
In particular $g^b(y')>y'$, otherwise the kneading sequence $\sigma^{b-j}(\hat{\alpha})$ would start with the word $\beta$. This contradicts the kneading information obtained from $f$. We proved that $g^b((0,y'))\supset (0,y')$ and hence the existence of a periodic point $q'\in (0,y')$ with $g^b(q')=q'$. 

\flushpar
Let $(x',0)$ be the maximal interval on which $g^a$ is monotone, the branch of type $\alpha$. In a similar way as above we show that $f^a((x',0))\supset (x',0)$ and the existence of a periodic point $p'\in (x',0)$ with $g^a(p')=p'$.

\bigskip

\flushpar
Left is to show that $(p',q',g^a,g^b)$ is an $(\alpha,\beta)$ renormalization.
First we will show that the orbits of $p'$ and $q'$ never enters $(p',q')$. 
 The kneading sequences of the periodic orbits of $p'$ and $q'$ are respectively $\alpha^{\infty}$ and $\beta^{\infty}$. The periodic orbits of $p$ and $q$ of $f$ do never enter in the interval $(p,q)$ and also have kneading sequences 
respectively $\alpha^\infty$ and $\beta^\infty$. This implies the following kneading information: for every $k\ge 0$ it is impossible that
$$
\alpha^{\infty}< \sigma^k(\beta^{\infty}), \sigma^k(\alpha^{\infty})< \beta^\infty.
$$ 
This kneading information implies that the orbits of $p'$ and $q'$ never enter $(p',q')$.

\flushpar
Left is to show that $g^b((0,q'))\subset [p',q']$ (and $g^a((p',0))\subset [p',q']$). Because $g$ has an $(\hat{\alpha},\hat{\beta})$ renormalization
there exists a branch $\hat{q}\in (0,z)\subset (0,q')$ with kneading sequence 
$\hat{\beta}=\beta\alpha\dots$. If $g^b(0_+)$ is left of $p'$ then also
$g^b((0,z)$ is left of $p'$ ( the orbit of $p'$ never enters $(p',q')$). According to the word $\hat{\beta}=\beta\alpha\alpha\dots$ we have to apply
a few times the branch $g^a|(x',0)$. Then $ g^{b+i\cdot a}((0,z))$ will still be on the left side of $p'$. This is impossible. Either $b+i\cdot a=\hat{b}$, in which case $g^{b+i\cdot a}(\hat{q})=\hat{q}$ but $g^{b+i\cdot a}(\hat{q})< p'$.
Or  $b+i\cdot a<\hat{b}$, in which case we have to apply the branch $g^b|(0,y)$
but  $g^{b+i\cdot a}((0,z))\cap (0,y)=\emptyset$.

\flushpar
We proved
that $(p',q',g^a,g^b)$ is a $(\alpha,\beta)$ renormalization.
\hfill\hfill\qed $\,\,$ (Proposition 3.3)
\enddemo

\vfill\eject
\bigskip
\centerline{\bf 4. Realization of Finite Combinatorics}
\bigskip

\proclaim { Proposition 4.1} 
     Let $V= [-1,0]\times [0,1]$ and $\Lambda \subset \bold R^2$ be homeomorphic to $V$. Let 
$$\Lambda \ni \lambda \mapsto f_\lambda \colon [-1,1] \to [-1,1]$$ be a family 
of Lorenz maps satisfying the following properties:
\item {1)} the branches $f_{\lambda,\pm}$ are $C^1$, $f_\lambda (\pm 1)= \pm 1$,  and there exists $K>0$ 
         such that $0<Df_\lambda (x)\leq K$ for all $x\neq 0$;
\item {2)} $Df_\lambda (x) \to 0$ as $x \to 0$; 
\item {3)} Let $F\colon \Lambda \to V$ be  defined by $F(\lambda)= (f_\lambda (0+),f_\lambda (0-))$ then $F(\partial \Lambda) \subset \partial V$  and 
the degree of the map $F\vert \partial \Lambda$ is different from zero. 

\flushpar
If $g\in \Cal{L}^r$, $r\ge 0$, is a simple Lorenz map with finite critical orbits then there exists $\lambda$ and a maximal semi conjugation $h$ from
 $f_\lambda$ to $g$. 
\endproclaim 

\flushpar
The proof of this Proposition is a corrected version of the proof 
of a similar statement for continous interval maps presented in [MS]. 

\demo { Proof} 
Let $g\in\Cal{L}^r$, $r\ge 0$ be a simple Lorenz map with finite critical orbits. 

\flushpar
Let $P(g)= \{z_1<\dots < z_l=0< z_{l+1}< \dots< z_k\}$ be the post-critical 
set of $g$ where  $-1 \leq z_1= g(0_+)\leq 0$ and $0 \leq z_k = g(0_-) \leq 1$ are
the critical values of $g$.  The order of these points and the mapping 
 $ c \colon \{ 1, 2, \dots ,l-1, l+1, 
\dots k\} \to \{ 1, \dots k\}$  defined by $g(z_i) = z_{c(i)}$ describes 
the combinatorics which we want to show can be realized by some map of 
our family. Consider the $k-1$ dimensional simplex $P= \{ x\in R^k \vert 
-1=x_0 \leq x_1 \leq \dots \leq x_l=0 \leq x_{l+1} \leq \dots \leq x_k \leq 
x_{k+1}=1\} $ and let $P_0$ be the interior of $P$.  
Let $pr \colon P \to V$ be the projection  $pr (x) = (x_1,x_k)$ and 
$$ Z = \{ (\lambda, x) \in \Lambda \times P \vert pr (x) = F(\lambda)\}$$ 
Let us consider the mapping $T\colon Z \to P$, $(\lambda, x) \mapsto y$, 
where $f_\lambda (y_i) = x _{c(i)}$  and $y_i$ has the same sign as $x_i$. 
It is clear that $T$ is well defined and continuous. To finish the proof 
we have to find $(\lambda,x) \in Z$ such that $T(\lambda,x)= x$. 

\flushpar
Let $\rho (x)= \text {min} \{ |x_i-x_{i+1}|; i=0, 1, \dots, k\}$ and $d(x,y)= 
\text {max} \{|x_i-y_i|; i= 1, \dots, k\} $. We need the 
following: 

\proclaim {Lemma 4.2} If $(\lambda_n, x(n))\in Z$ is such that $\rho (x(n)) \to 0$ 
then 
$$ \lim_{n\to \infty} \frac {d(T(\lambda_n,x(n)), x(n))}{ \rho (x(n))} = \infty $$
\endproclaim 

\demo {Proof} 
Suppose, by contradiction, that there exists $K_0>0$ and a sequence $(\lambda_n, x(n)) \in Z$ such that $\rho(x(n)) \to 0 $ and
 $$
d(x(n),y(n)) \leq K_0 \rho (x(n))
$$ 
where $y(n)= T(\lambda_n, x(n))$. 

\proclaim { Claim 1} There exists constant $ K_1>0$ such that if $x_i(n), x_j(n)$ have the same sign then 
$$ |x_{c(i)} (n) - x_{c(j)}(n)| \leq  K|x_i(n)-x_j(n)| + K_1 \rho (x(n)).$$
\endproclaim

\flushpar
Since $Df_\lambda (t) \leq K$, it follows from the Mean Value Theorem that 
$$
\aligned
|x_{c(i)}(n)-x_{c(j)}(n)| 
&\leq K |y_i(n)-y_j(n)| \\
&\leq K|x_i(n)-x_j(n)|+K|y_i(n)-x_i(n)| + K|y_j(n)-x_j(n)|\\
&\leq K|x_i(n)-x_j(n)| + 2K_0\rho(x(n))
\endaligned
$$ 
and the claim is proved. 

\bigskip

\flushpar
From Claim 1 we get by induction, for all $s$ there exists constant $C_s>0$ such that 
 if $x_{c^k(i)}(n)$ and $x_{c^k(j)}(n)$ have the same sign for $k<s$ then, 
$$ 
|x_{c^s(i)}(n)- x_{c^s(j)}(n)| \leq C_s |x_i(n)- x_j(n)| + C_s \rho (x(n)) $$ 

\proclaim{Claim 2} 
There exists $s_0$ such that for all $m\in \{1, \dots, k\}$ there exists $s\leq 
s_0$ with $0 \in [x_{c^s(m)}(n), x_{c^s(m+1)}(n)] $. 
\endproclaim

\flushpar
Consider the set $\Cal P$ of pairs $(z_i, z_j)$ such that $0\not\in [z_i,z_j]$. 
If the claim is false then, since $\Cal P$ is finite, there exists a periodic 
pair in $\Cal P$. This implies the existence of a periodic homterval non-essential for $g$.

\flushpar
To finish the proof of the lemma let $m$ be such that $\rho (x(n))= |x_m(n)-x_{m+1}(n)|$ ( $m$ depends on $n$). Let $s\leq s_0$ be such that 
$0 \in [x_{c^s(m)}(n), x_{c^s(m+1)}(n)]$ and $x_{c^k(m)}(n), x_{c^k(m+1)}(n)$
have the same sign for
$k<s$. 

\flushpar
From Claim 2 we get some $t$, depending on $n$ such that $|x_t(n)-x_{t+1}(n)| \leq C_s \rho (x(n))$ and either $x_t(n)=0$ or $x_{t+1}(n)=0$. 
Since $|f_{\lambda_n}([y_t(n), y_{t+1}(n)])| = |[x_{c(t)}(n), x_{c(t+1)}(n)]| \geq \rho(x(n))$ we have that 
$Df_{\lambda _n} (\theta_n) |y_t(n)- y_{t+1}(n)| \geq \rho (x(n))$ for some 
$\theta_n \in [y_t(n),y_{t+1}(n)]$. 

\flushpar
Consider the case when $x_{t+1}(n)=0$. Then also $y_{t+1}(n)=0$. Therefore, 
$$
\aligned
d(y(n), x(n)) &\geq |y_t(n)-x_t(n)|\\
              &\geq \frac {1}{Df_{\lambda_n}(\theta_n)}\cdot \rho(x(n))
                    -C\cdot \rho (x(n))\\
              &\geq (\frac {1}{Df_{\lambda_n}(\theta_n)} -C)\rho (x(n))
\endaligned
$$ 
On the other hand, $|y_j(n)-x_j(n)|\leq d(y(n),x(n)) \leq K_0 \rho(x(n))\to 0$ as $n\to \infty$. This implies that $Df_{\lambda_n}(\theta _n) \to 0$. 
This contradicts $d(x(n),y(n))\le K_0\rho(x(n))$. 

\hfill\hfill\qed $\,\,$ (Lemma 4.2)
\enddemo

\flushpar
 Assume that there is no fixed point for the mapping 
$T\colon Z \to P$.

\proclaim { Lemma 4.3} There exists a continuous one parameter family of 
continuous mappings $G_t \colon Z \to P$ such that 
\item 1) $G_0 (\lambda, x) = x$
\item 2) $G_t (\lambda, x) = x$ if $x \in \partial P$ 
\item 3) $G_1 (Z) \subset \partial P$. 
\endproclaim 

\demo {Proof}  Let $D=\{x\in \Bbb{R}^{k-1}| \sum_{i=1}^{k-1} x_i^2\le 1\}$ be the closed unit ball in $\bold R^{k-1}$
and let $\Phi\colon P \to D$ be a bi-Lipschitz homeomorphism. 
Let $\tilde T \colon Z \to D$ be the mapping $\tilde T= \Phi\circ T$. Since 
$\Phi $ is bi-Lipschitz and the distance from $x$ to the boundary of $P$ is 
equal to $\frac12\rho(x)$ 
we have that, if $(\lambda_n,x(n))\in Z$ is such that 
$\rho (x(n)) \to 0$ then 
$$ 
\frac { d( \Phi (x_n), \tilde T (\lambda_n, x(n))} {d(\Phi (x(n)), \partial D)} \to \infty 
$$ 

\flushpar
Now define $\tilde G_1 (\lambda, x)$ as the point in the boundary of $D$ 
such that $\Phi(x)$ lies in the line segment bounded by $\tilde T (\lambda, x)$ 
and $\tilde G_1(\lambda, x)$. It follows that $\tilde G_1$ is continuous and
$\tilde {G}_1(\lambda,x)=\Phi(x)$ for $x\in \partial P$.
Define $\tilde G_t= t \tilde G_1 + (1-t) \Phi$ and $G_t = \Phi^{-1} \circ 
\tilde G_t$. This clearly satisfies the conditions of the Lemma. 
\hfill\hfill\qed $\,\,$ (Lemma 4.3)
\enddemo

\proclaim {Lemma 4.4}  Let $h\colon V \to \Lambda$ be a homeomorphism such that $\pi_i(x)=0 \Rightarrow \pi_i (F(h(x)) =0$ for $i=1, 2$, where $\pi_i\colon V 
\to \bold R$ is the projection in the i-th coordinate. Then there exists a 
mapping $H\colon P \to Z$ with the following properties:
\item {1)} $H$ is continuous and the restriction of $H$ to the interior of 
$P$ is a homeomorphism;
\item {2)} $\pi \circ H= h\circ pr $ where $\pi \colon Z \to \Lambda$ is the 
projection $(\lambda,x) \mapsto \lambda$; 
\item {3)} the restriction of $H$ to each fiber $pr^{-1}(v)\cup int(P)$ is a 
diffeomorphism onto the fiber $\pi^{-1} (h(v))\cup \text{int}(Z)$. 
\endproclaim 

\demo {Proof} Let $v=(v_1,v_2)$ and $F(h(v))=(\tilde{v_1},\tilde{v_2})$. Let $H_v:\Bbb{R}\to\Bbb{R}$ be the Moebius transformation that maps $v_1$ to $\tilde{v_1}$, $0$ to $0$ and $v_2$ to $\tilde{v_2}$. For $x\in pr^{-1}(v)\cap\text{int}(P)$ define $H(x)=(h(v),(H_v(x_1),\cdots, H_v(x_k)))$. 
\hfill\hfill\qed $\,\,$ (Lemma 4.4)
\enddemo

\flushpar
Choose a homeomorphsim $h:V\to \Lambda$ with the property needed to apply Lemma 4.4 and let $H:P\to Z$ be the corresponding map from Lemma 4.4. Let $\hat G_t= G_t\circ H \colon P \to P$. 
We have that 
\item {1)} $\hat G_t$ is continuous and depends continuously on $t$;
\item {2)} $\hat G_t ( \partial P )\subset \partial P$ and $\hat G_t ( \text{int }(P)) \subset \text{int}(P)$; 
\item {3)}   $\hat{G}_t(x)=G_0(H(x))$ for all $x\in \partial P$;
\item {4)} $\hat G_1( P) \subset \partial P$ ; 
\item {5)} The degree of $\hat G_0$ is equal to the degree of $F$. 

\flushpar
The only statement that needs a proof is 5). Notice that
$$
pr\circ \hat{G}_0=(F\circ h)\circ pr
$$
and the restriction of $\hat{G}_0$ to $pr^{-1}(v)\cap \text{int}(P)$ is a diffeomorphism onto $pr^{-1}(F(h(v))\cap \text{int}(P)$. Hence the degree of
$\hat{G}_0$ is equal to the degree of $F\circ h$. Which is equal to the degree of $F$. Here we are using the following topological fact. Let $\Phi:D\to D$ be a continuous mapping that maps $\partial D$ onto $\partial D$, $\text{int}(D)$ onto $\text{int}(D)$ and $\Phi |_{\text{int}(D)}$ is smooth. Then the degree of $\Phi |_{\partial D}$ equals the degree of $\Phi|{\text{int}(D)}$, (see [D, p 67]). This is a contradiction.

\bigskip

\flushpar
It follows that there is a map $f_\lambda$ in the family which has also periodic critical orbits. Moreover, the combinatorics of these critical orbits are the same as the combinatorics of the critical orbits of $g$.

\flushpar
Let $h_n$, $n\ge 0$ be the homeomorphism that maps $B_n(f)$ into $B_n(g)$ 
given by 
Lemma 2.1. Because $g$ is a simple Lorenz map it follows that $h_n$ converges to a maximal semi conjugacy $h$ from $f$ to $g$. 
\hfill\hfill\qed $\,\,$ (Proposition 4.1)
\enddemo

\bigskip
\centerline{\bf 5. Archipelagoes in the parameter plane}
\bigskip

\flushpar
In this section we will study the parameter plane of a given 
monotone Lorenz family $F:[0,1]\times [0,1]\to \Cal{L}^r$, $r\ge 0$. The main object of our
study is to understand the topological structure of Archipelagoes. 

\flushpar  
Let us concentrate on the archipelago $A=A_{\alpha,\beta}$. For every
$\lambda\in A$ the  domain of an $(\alpha,\beta)-$renormalization is denoted by $(p_\lambda,q_\lambda)$. 
Furthermore let $a=|\alpha|$ and $b=|\beta|$.

\proclaim{Proposition 5.1 (Island Structure)} 
Let $I$ be an island of the archipelago $A$. Then
\parindent=15pt
\item{1)} For all $u\in U$, $M_u\cap A_{\alpha,\beta}=M_u\cap I$.
\item{2)} There exists an interval $[u_1,u_2]\subset U$ and Lipschitz
functions $\partial_+,\partial_-:[u_1,u_2]\to M$ with Lipschitz constant $1$
and $\partial_-(x)<\partial_+(x)$, $x\in (u_1,u_2)$, such that
$$
I=\{(u,m)| u\in (u_1,u_2) \text{ and } \partial_-(u)<m<\partial_+(u)\}.
$$
Furthermore $\partial_-(u_1)=\partial_+(u_1)$ and $\partial_-(u_2)=\partial_+(u_2)$.
\item{3)} If $\lambda=(u_1,\partial_-(u_1))$ or $\lambda=(u_2,\partial_-(u_2))$, one 
of the
{\it extremal} points of the island then
$$
f_\lambda^a(0_-)=0 \text{ and } f_\lambda^b(0_+)=0
$$
or
$$
f_\lambda^a(0_-)=q_\lambda \text{ and } f_\lambda^b(0_+)=p_\lambda.
$$
The first possibility is called a {\it trivial extremal point}, the second a
{\it full-branch extremal point}.
\item{4)} For all $\lambda\in \partial_+$
$$
f_\lambda^a(0_-)=q_\lambda \text{ or } f_\lambda^b(0_+)=0.
$$
\item{5)} For all $\lambda\in \partial_-$
$$
f_\lambda^a(0_-)=0 \text{ or } f_\lambda^b(0_+)=p_\lambda.
$$
\endproclaim

\demo{proof} Let $(p_m,q_m,f_m^a,f^b_m)$ be the renormalization of $f_m\in I\cap M_u$. Observe that $m\mapsto q_m$ and $m\mapsto p_m$ are strictly monotone decreasing and $m\mapsto f_m^a(0_-)$ and $m\mapsto f_m^b(0_+)$ are strictly monotone increasing. This implies that once we arrive at the boundary of $I\cap M_u$, by increasing $m$ up to $m_+$ we are in the situation 
$$
f^a_m(0_-)=q_m \text{ or } f^b_m(0_+)=0.
$$
In particular for any $m>m_+$ we have 
$$
\theta^-(m')>\alpha\beta^\infty \text{ or } \theta^+(m')\ge \beta^\infty,
$$
and hence $m'\notin A$. Similarly we study what happens when $m$ decreases up to the boundary of $I\cap M_u$. We proved that 
$$
A\cap M_u= I\cap M_u
$$
and that this intersection is an interval $(\partial_-(u),\partial_+(u)$. 
The boundary of $I\cap M_u$ consists of points $\partial_-(u)$ and points $\partial_+(u)$: the boundary of $I$ consists of two parts, $\partial_- I$ and $\partial_+ I$. Moreover
$$
f^b_{\partial_-(u)}(0_+)=p_{\partial_-(u)}  \text{ or } 
f^a_{\partial_-(u)}(0_-)=0
$$
and 
$$
f^a_{\partial_+(u)}(0_-)=q_{\partial_+(u)} \text{ or }
f^b_{\partial_+(u)}(0_+)=0.
$$
The island $I$ is open and connected. Hence the set of values
$u$ for which $M_u\cap I\ne \emptyset$ is an open interval, say $(u_1,u_2)$.
Now observe that for any $z\in \partial_\pm I$ we have 
$$
\partial_\pm I\subset B_z,
$$
where $B_z$ is the complement of the cones $C_z^+$, $C_z^-$ defined in section 2.
In particular the functions $\partial_\pm:(u_1,u_2)\to M$ are Lipschitz.

\flushpar
Left is to explain the behavior of the boundary points $\partial_\pm(u)$ when
$u$ tends to $u_1$ or $u_2$. For every $u\in (u_1,u_2)$ we have $\partial_-(u)<\partial_+(u)$. Let us assume that also $\partial_-(u_1)<\partial_+(u_1)$. The archipelago $A$ is closed. Hence for any 
$(u_1,m)$ with $\partial_-(u)<m<\partial_+(u)$ there is renormalization
$(p_m,q_m,f^a_m, f^b_m)$. Because $(u_1,m)\in C_-((u_1,\partial_+(u_1))$ and $(u_1,m)\in C_+((u_1,\partial_-(u_1))$ we have $f^a_m(0_-)< q_m$ and $f^b(0_+)>p_m$. Because of Lemma 3.1 the point $(u_1,m)$ is in the interior of $I$. This is a contradiction, proving that $\partial_-(u_1)=\partial_+(u_1)$.
\hfill\hfill\qed $\,\,$ (Proposition 5.1)
\enddemo

\flushpar
We have to distinguish special points on the boundary of islands. One
type of special points are the extremal points, discussed in the above 
proposition. The other type of special points are vertices. Let $\partial_-$
be the lower part of the boundary of the island $I$.
A point $\lambda\in \partial_-$ in the lower boundary  is called a 
{\it vertex} if
$$
f_\lambda^a(0_-)=q_\lambda \text{ and } f_\lambda^b(0_+)=0.
$$
A vertex in the upper boundary is defined similarly.

\proclaim{Lemma 5.2} Every archipelago contains only finitely many islands
with a vertex.
\endproclaim

\demo{proof}
Suppose the archipelago $A_{\alpha,\beta}$ has infinitely many island $I_n$ which
has a vertex $v_n\in \partial_+^n$. Here $\partial_+^n$ denotes the upper 
boundary of $I_n$. We may assume that $v_n\to v$. Since archipelagoes are closed sets, 
$v\in A_{\alpha,\beta}$. 

\flushpar
The topological type of the
maps in the vertices are all the same and it follows easily that the map corresponding to the parameter $v$ inherits this type also. So
$$
f^a_v(0_-)=q_v \text{ and } f^b_v(0_+)=0.
$$
By decreasing both parameters a little bit, we observe that $v$ is actually
a vertex of some island $I$ in the same archipelago $A_{\alpha,\beta}$.

\flushpar
From Proposition 5.1 we get a neighborhood $(a,b)\times M$ of $I$ such that
the archipelago intersects this neighborhood only in $I$, there are no other 
islands intersecting this neighborhood. But the vertex $v$ of $I$ lies in 
this neighborhood, it cannot be accumulated by the islands $I_n$. 
This is a contradiction.
\hfill\hfill\qed $\,\,$ (Lemma 5.2)
\enddemo

\proclaim{Lemma 5.3} Let $\lambda$ be a parameter for which $0_-$ is periodic.
Say with period $a=|\alpha|$, where $\alpha$ is the word describing the 
combinatorics of the orbit of  $0_-$. Denote the set of parameters for 
which $0_-$ is periodic with the same combinatorics $\alpha$ by $\alpha_-$.

\flushpar
Then there exists an interval $I\subset U$ and 
a Lipschitz function $\gamma: I\to M$ with Lipschitz constant $1$ such that
$I\times M$ is a neighborhood of $\lambda$ and
$$
graph(\gamma)\cap (I\times M)= \alpha_-\cap (I\times M).
$$
\endproclaim

\demo{Proof} Let $z=(u,m)\in\alpha_-$. The monotonicity of the family implies
$$
\alpha_-\subset B_z
$$
and
$$
\aligned
K^-(z')&< \alpha, \text{    } z'\in \partial C^-_z\setminus \{z\}\\
K^-(z')&> \alpha, \text{    } z'\in \partial C^+_z\setminus \{z\}.
\endaligned
$$
For a small enough nieghborhood $I\times N$ of $z\in \alpha_-$ the orbit of $0_-$ will have the combinatorics given by $\alpha$ except for the last symbol. For each $u'\in I$ let $m'_\pm$ be such that $(u',m'_\pm)\in \partial C^\pm(z)$. According to the above observation we have 
$$
\aligned
K^-((u',m'_-))&< \alpha\\
K^-((u',m'_+))&> \alpha.
\endaligned
$$
In particular for each $u'\in I$ there exists a unique $m'(u')\in (m'_-,m'_+)
\subset M_{u'}$
with
$$
(u',m'(u'))\in \alpha_-
$$ 
the set $\alpha_-\cap I\times M$ is the graph of the function $m':I\to M$.
The first property discussed in this proof, namely $\alpha_-\subset B_z$
for each $z\in \alpha_-$ implies that this function $m'$ is $1-$Lipschitz.
\hfill\hfill\qed $\,\,$ (Lemma 5.3)
\enddemo

\proclaim{Lemma 5.4} Let $\lambda$ be a parameter for which the map 
$f_\lambda$ has two periodic points $p<0<q$. Say with period $a=|\alpha|$ and 
$b=\beta$ respectively.
Furthermore assume
\parindent=15pt
\item{1)} $(p,q)\cap(orb(p)\cup orb(q))=\emptyset$.
\item{2)} $f^a|[p,0]$ and $f^b|[0,q]$ are monotone.
\item{3)} $Df^a(p)>1$ and $Df^b(q)>1$. 
\item{4)} $f^a([p,0])=[p,q]$

\flushpar
The word $\alpha\beta^\infty$ describes the combinatorics of the orbit of 
$0_-$. Denote the set of parameters for which the corresponding map has the 
properties 1), 2), 3) and 4)  by $\alpha\beta^\infty_-$.

\flushpar
Then there exists an interval $I\subset U$ and 
a Lipschitz function $\gamma: I\to M$ with Lipschitz constant $1$ such that
$I\times M$ is a neighborhood of $\lambda$ and
$$
graph(\gamma)\cap (I\times M)= \alpha\beta^\infty_-\cap (I\times M)
$$
\endproclaim

\demo{Proof} Let $z\in \alpha\beta^\infty_-$. The properties 1), 2) and 3)
hold in a neighborhood of $z$. The proof continues in this neighborhood as the 
proof of Lemma 5.3.
\hfill\hfill\qed $\,\,$ (Lemma 5.4)
\enddemo

\flushpar
Similar statements hold for the combinatorics of the orbit of $0_+$. The
corresponding sets of parameters will be denoted by resp. $\beta_+$ and 
$\beta\alpha^\infty_+$.

\flushpar
The uniqueness part of the above Lemma allow us to consider maximal arcs: the 
connected components of the sets $\alpha_-$ and $\alpha\beta^\infty_-$ are
graphs of Lipschitz functions.

\proclaim{Lemma 5.5} Let $\gamma: (t_0,t_1)\to M$ be the Lipschitz function 
whose graph is a component of $\alpha\beta^\infty_-$. Then
\parindent=15pt
\item{1)} This function can be extended to a Lipschitz function on $[t_0,t_1]$.
\item{2)} There exists a function $\gamma_2:[t_0,t_1]\to M$ whose graph is 
contained in $\beta_+$.
\item{3)} If there exits $(u,m)\in (graph(\gamma)\cap A_{\alpha,\beta})\setminus \beta_+$ then 
$\beta_+$ intersects each component of $graph(\gamma)\setminus\{\lambda\}$ 
unless such a component terminates in the boundary of the parameter domain. 
Moreover $\gamma_2(u)>\gamma(u)$, the $\beta_+$ curve lies above the point
$(u,m)$ in the $\alpha\beta^\infty_-$ curve.

\flushpar
In particular a component of $\alpha\beta^\infty_-\setminus \{(u,m)\}$ 
which does not terminate in the boundary of the parameter domain contains a 
vertex in the upper boundary of an island of the archipelago $A_{\alpha,\beta}$.
\endproclaim

\demo{Proof} Remember that the arc $\gamma:(t_1,t_2)\to M$ satisfies
$$
\text{graph}(\gamma)\subset B_{(t,\gamma(t))}
$$  
for all $t\in (t_1,t_2)$, $\gamma$ is $1-$Lipschitz. Observe that the length of the interval $T_t=M_{t_1}\cap B_{(t,\gamma(t))}$ tends to zero, when $t\to t_1$.
Moreover these intervals $T_t$ with $t>t_1$ are nested: $T_s\subset T_t$ whenever $t_1<s<t$. The intersection of the intervals $T_t$ defines the continuous extension of the arc $\gamma$.  

\flushpar
To proof 2) it is enough to show that $\beta_+\cap M_t\ne \emptyset$ for each $t\in (t_1,t_2)$. Let $t\in (t_1,t_2)$ and $\{\lambda\}= \alpha\beta_-^\infty\cap M_t$. 

\demo{case I: $f^b(0_+)>0$} Consider the maps in $M_t$ with the properties
\parindent=15pt
\item{1)} there exists a branch $T=(0,x)$ of type $\beta$,
\item{2)} there exists $x_0\in T$ with $f^b(x_0)<x_0
$,
\item{3)} $f^b(0_+)>0$.

\flushpar
Observe that property 2) and 3) imply the existence of a periodic attractor in $T$ which attracts the orbit of $0_+$. In particular, the orbit of $0_+$ is infinite. This implies that $0\notin \partial_- f^i(T)$ for each $i\le b$.
Consequently, if $(t,m)$ is a point with this properties there is an small
interval $\{t\}\times (m-\epsilon,m]$ of maps with the three property.

\flushpar
Observe that the map in $\lambda$ has the three properties. Let $H=\{t\}\times (m_1,m_0]$ be the maximal interval of maps with the above properties. We are going to show that in $(t,m_1)$ the map is in $\beta_+$.

\flushpar
Observe that by decreasing $m\in (m_1,m_0]$ we see that the interval $T$ is increasing in length. Because the family is monotone, the same point $x_0$ in  $T$ which was moved to the left persists to be in $T$ and will be moved to the left. As before we see that the orbit of the interval $T$ will
never hit $0_+$.

\flushpar
We showed that property 1) and 2) hold also in $(t,m_1)$. Hence, in the boundary of $H$ property 3) has to be violated. Otherwise we could decrease $m$ slightly more. 
\enddemo

\demo{case II: $f^b(0_+)<0$} Consider the maps in $M_t$ with the properties
\item{1)} there exists a branch $T=(0,x)$ of type $\beta$,
\item{2)} $f^b(x)>x$,
\item{3)} $f^b(0_+)<0$.

\flushpar
Properties 2) and 3) imply that there is a hyperbolic periodic point $q\in T$
of type $\beta$. Consequently, the branch of type $\beta$ will persist under small perturbations. The collection of above maps is open.

\flushpar
Observe that the map in $\lambda$ satisfies the above properties. We have $f^b(x)>x$
otherwise the orbit of $0_-$ could not be trapped in the orbit of type $\beta$,
which is the case in $\lambda\in \alpha\beta^\infty_-$. Let $H=\{t\}\times [m_0,m_1)$ be the maximal interval of maps with the above properties. We are going to prove that $(t,m_1)\in \beta_+$.

\flushpar
Observe that by increasing $m\in [m_0,m_1)$ the branch $T$ will be decreasing, a consequence of the fact that the family is monotone. However the periodic orbit
persists and there has to be a periodic orbit at $(t,m_1)$. Maybe not hyperbolic anymore. In particular, there is a branch $T=(0,x)$ of type $\beta$ in $(t,m_1)$.

\flushpar
Again from the monotonicity of the family we get that $f^b(x)>x$ in $(t,m_1)$.
The only way to reach the boundary of $H$ is to violate property 3): $(t,m_1)\in \beta_+$. 
\enddemo

\bigskip

\flushpar
To prove 3) Let $(t,m)=\lambda\in \text{graph}(\gamma)\cap A_{\alpha,\beta})\setminus \beta_+$. The $M_t$ contains a point in $\beta_+$. 
Because $(t,m)\in A_{\alpha,\beta}\setminus \beta_+$, there is a renormalization, we get that $\gamma_2(t)>\gamma(t)$. To finish the proof of $3)$ we have to show that $f^b(0_+)>0$ in $(t_1,\gamma(t_1)$.

\flushpar
Because $(t_1,\gamma(t_1))$ is a boundary point we get that $Df^b(q)=1$. Otherwise we could extend $\alpha\beta_-^\infty$. The map under consideration has negative Schwarzian derivative. Hence, the neutral periodic point attracts a critical point with an infinite orbit. The orbit of $0_-$ is trapped in the orbit of $q$. Only the orbit of $0_+$ can be attracted towards $q$. In particular $f^b(0_+)>0$.
\hfill\hfill\qed $\,\,$ (Lemma 5.5)
\enddemo

\flushpar
A similar Lemma holds for the combinatorics of $0_+$. In the next Lemma we will see
that the symmetry breaks down. Moving up or down has well understood 
consequences on the combinatorics of the map. On the other hand moving one 
branch up and the other down are the deformations which are difficult to
understand, the deformations in directions parallel to U. The next Lemma
indicates a difference between such movements to the right and to the left. 

\proclaim{Lemma 5.6} Let $\lambda\in 
\alpha\beta^\infty_-\cap \beta\alpha^\infty_+$, say $\lambda=(u,m)$. Then 
there exists a parameter $\mu=(t,m')$ such that
\parindent=15pt
\item{1)} $t\le u$.
\item{2)} $\mu$ is the right extremal point of an island $I$ of the 
archipelago $A_{\alpha,\beta}$. Moreover it is a full branch extremal point.

\item{3)} The island $I$ has a vertex in its upper (and lower) boundary.
\endproclaim

\bigskip
\centerline{\psfig{figure=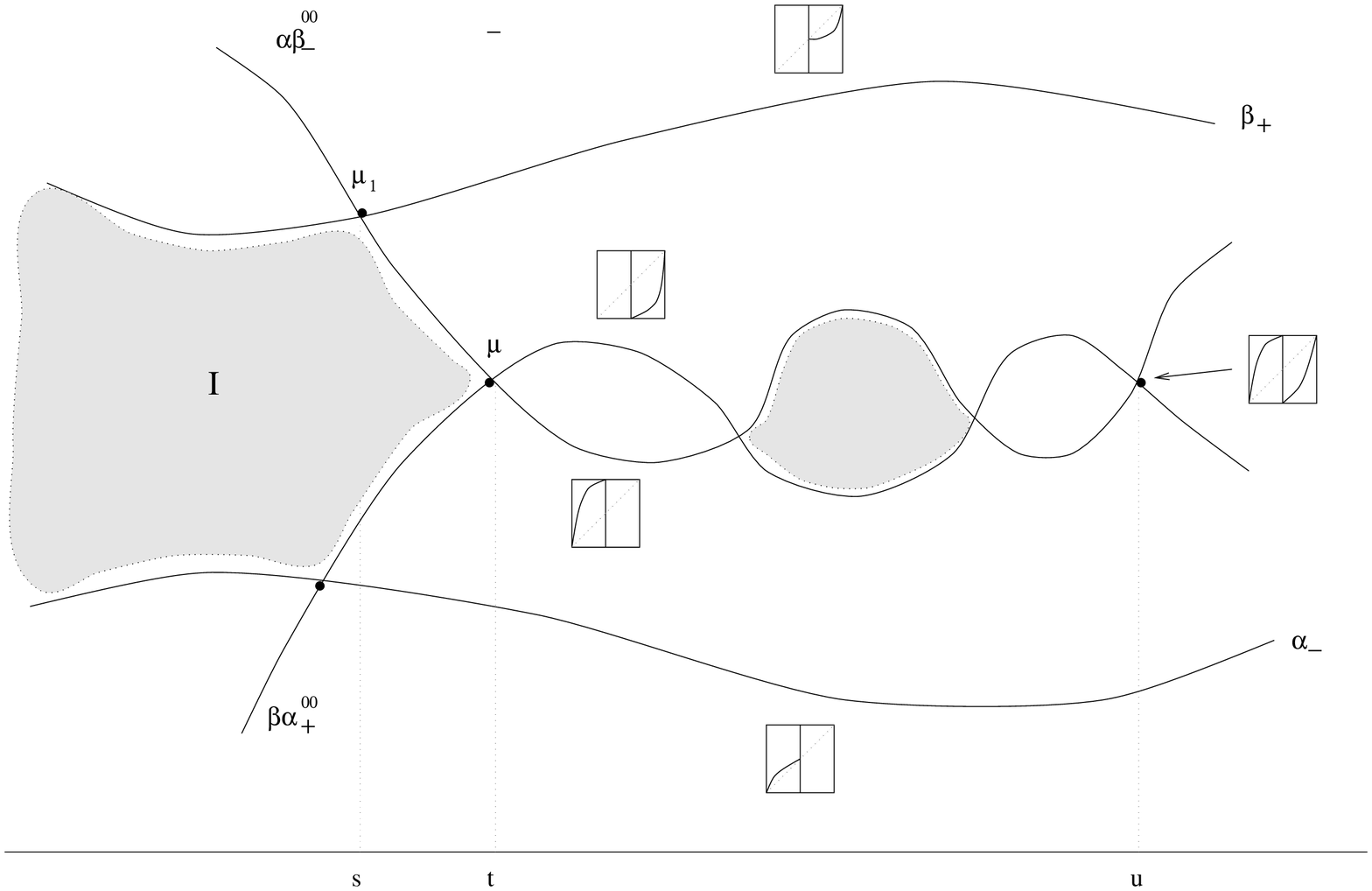,width=10cm}}
\smallskip
\centerline{\bf Figure 2 Illustration to the proof of Lemma 5.6}
\bigskip

\demo{proof} Let $\gamma_+:(t_+,u]\to M$ be the function whose graph is the 
component of $\alpha\beta^\infty_-\setminus \{\lambda\}$ which touches 
$\lambda$ from the left side. Similarly $\gamma_-:(t_-,u]\to M$ be the
function whose graph is the component of 
$\beta\alpha^\infty_+\setminus \{\lambda\}$ which also touches $\lambda$ 
from the left side. 

\flushpar
Let $t\in [t_+,u]\cap [t_-,u]$ be minimal such that $\gamma_+(t)=\gamma_-(t)$,
describing the left most intersection point of $\alpha\beta^\infty_-$ with
$\beta\alpha^\infty_+$. Let $\mu$ be the intersection point. Observe that
Lemma 5.4 can be applied in this point, giving us local extensions of the arcs
$\alpha\beta^\infty_-$ and $\beta\alpha^\infty_+$. So $\mu$ is a point in the 
interior of both arcs.

\flushpar
Actually $\mu\in A_{\alpha,\beta}$, the map is of full branch type. So we can 
apply Lemma 5.5: let $s\le t$ be maximal such that 
$\mu_1=(s,\gamma_+(s))\in \beta_+\cap\alpha\beta^\infty_-$. Observe that $\alpha\beta_-^\infty$ can not terminate to the left in the boundary of 
parameter space.

\flushpar
Clearly $\mu_1$
is a vertex point in the upper boundary of some island 
$I\subset A_{\alpha,\beta}$. Observe that the island $I$ intersects 
$[s,t]\times M$ only below the graph of $\gamma_+|[s,t]$, moving up a point
in this graph would destroy the renormalization immediately.

\flushpar
The aim is to show that the graph of $\gamma_+|[s,t]$ is part of the upper
boundary of the island $I$. From Lemma 5.5 we get that the arc $\beta_+$ lies 
above the point $\mu$. Furthermore $\mu_1$ is the right most intersection of 
the graph of $\gamma_+|(t_+,t]$ with the $\beta_+$ arc. So the $\beta_+$ arc 
in the strip $(s,t)\times M$ lies above the graph of $\gamma_+$, which lies
above the island $I$. So the boundary points of the island $I$ in the strip 
$(s,t)\times M$ lie in the graph of $\gamma_+$, boundary points in the upper
boundary can only be of two types $\alpha\beta^\infty$ or $\beta_+$.

\flushpar
Left is to show that the boundary actually extends up to $\mu$. Clearly
the piece of the boundary in the strip can not have vertices anymore, for 
this you need points in $\beta_+$. So the graph contains the right extremal 
point of the island. This extremal point can not be trivial, for this you need
a point in $\beta_+$. So it is of full branch type. The only parameter of this
type in the graph of $\gamma_+|[s,t]$ is $\mu$. Clearly the island does not 
extend beyond $\mu$, the archipelago intersects $\{t\}\times M$ only in $\mu$.
\hfill\hfill\qed $\,\,$ (Lemma 5.6)
\enddemo

\bigskip
\centerline{\bf 6. Proof of Full-Island Theorem 1.9}
\bigskip

\flushpar
Fix a monotone Lorenz family $F:[0,1]\times [0,1]\to \Cal{L}^r$, $r\ge 3$. Let $A=A_{\alpha,\beta}\subset [0,1]\times [0,1]$ be an the $\alpha,\beta-$archipelago.

\proclaim{Proposition 6.1} The archipelago $A$ has an island $I\subset A$ which has a trivial and a full branch extremal point.
\endproclaim 

\demo{proof} By Proposition 4.1 there exists $\lambda\in 
\alpha\beta^\infty_-\cap \beta\alpha^\infty_+$. Now Lemma 5.6 gives an
island which has a vertex and its right extremal point is of full branch 
type. Because there
are only finitely many islands with a vertex we can take the left most
island $I\subset A_{\alpha,\beta}$ which has a vertex and whose right 
extremal point is of full branch type.

\flushpar
Left is to show that the left extremal point of this island is trivial.
Suppose not and apply Lemma 5.6 again: we will find an island left of $I$
which has a vertex and whose right extremal point is of full branch type.
But $I$ was the left most island with these properties. This is a 
contradiction.
\hfill\hfill\qed $\,\,$ (Proposition 6.1)
\enddemo

\proclaim{Lemma 6.2} For each Lorenz map $f\in \Cal{L}^2$ there exists a simple Lorenz map
$\hat f$ and a maximal semi conjugation from $f$ to $\hat{f}$.
\endproclaim

\demo {Proof} Define $x\sim y$ if the there exists a countable closed set
$C\subset [x,y]$ such that each connected component of $[x,y]\setminus C$ is a homterval.
Clearly, it follows that $\sim$ is an equivalence relation and the 
equivalence classes are either points or closed intervals. 
Also $f$ maps equivalence classes into equivalence classes: if $[x]$ is the equivalence class of $x$ and $0\notin \text{\it int}[x]$ then
$$
f([x])\subset [f(x)].
$$
Therefore, the quotient space 
$[P,Q]/ \sim$ is an interval and $f$ induces a Lorenz map on this 
interval which does not have homtervals. 

\flushpar
The semi conjugacy being the quotient map $h$. Left is to show that $h$
is a maximal semi conjugacy. Suppose there exists $y\in [-1,1]$ such that
$h^{-1}(y)$ is not a point neither a homterval.  In this case $T=h^{-1}(y)$
contains a closed countable set $C$ such that the connected components of $T\setminus C$ are all homtervals. Because $C$ is countable and closed it has an
isolated point $c\in C$. In particular there are two homtervals 
$I_1, I_2\subset T$ with $c$ as common boundary point. Clearly, this common boundary point is a preimage of
$0$, say $f^n(c)=0$. In particular $f$ has two critical homtervals $f^n(I_1)$ and $f^n(I_2)$. 

\flushpar
A consequence of these two critical homtervals and the fact that $f$ is $C^2$
is that $f$ does not have wandering intervals. If $f$ would have a wandering 
interval then we could modify the map on
$f^n(I_1)\cup f^n(I_2)$ to obtain a smooth bimodal map with a wandering interval.
In [MMS] this is proved to be impossible.

\flushpar
Every homterval of $f$ eventually falls into a periodic homterval. We proved that $h$ is a maximal semi conjugacy.
\hfill\hfill\qed $\,\,$ (Lemma 6.2)
\enddemo

\parindent=20pt
\flushpar
Let us sumarize the possibilities in the above Lemma. 
\item{1)} If $f$ has at most one 
critical homterval then $h^{-1}(y)$ is a point or a homterval for each $y\in [-1,1]$.
\item{2)} If $f$ has two critical homtervals, say $L$ and $R$ then
\item{2a)} The two intervals $L$ and $R$ are periodic with distinct orbits or
\item{2b)} The two intervals are periodic with the same orbit or
\item{2c)} One interval is periodic and the other is eventually periodic or
\item{2d)} The two intervals are eventually periodic.

\flushpar
The notion of maximal semi conjugacy is needed to collapse al these periodic and
eventually periodic homtervals. If the map $f$ would have at most one critical homterval the usual equivalence relation $x\sim y$ if an only iff $[x,y]$
is a homterval would give the maximal semi conjugacy to a simple Lorenz map.

\proclaim { Theorem 1.9} Let $I\subset A_{\alpha,\beta}$ be an island which has a trivial and a full branch extremal point. The Lorenz Family $G:I\to \Cal{L}^r$
defined by
$$
I\ni \lambda \mapsto (p_\lambda,q_\lambda, f_\lambda^a,f_\lambda^b)\in 
\Cal{L}^r
$$
is a full family.
\endproclaim 

\demo{Proof} Let $f$ be a Lorenz map. We would like to find a map $g=G(\lambda)$
which is essentially conjugated to $f$. By Lemma 6.2 we get a 
simple Lorenz map $\hat{f}$ and a maximal semi conjugation $h_1$ 
from $f$ to $\hat{f}$. 

\flushpar
Let $K^\pm_n=K^\pm_n(\hat{f})$, $n\ge 0$. Proposition 4.1 can be applied to the 
family $G$. Hence for each $n\ge 0$ there exists $\lambda_n\in I$ such that
$K^\pm_n=K^\pm_n(g_n)$, where $g_n=G(\lambda_n)$. By Lemma 2.1 we get homeomorphisms $h_n:[-1,1]\to [-1,1]$ preserving the combinatorics up to time $n$:
$$
h_n: \Bbb{B}_n(g_n)\to \Bbb{B}_n(\hat{f}).
$$
We may assume that a subsequnce of $g_n$ converges to $g=G(\lambda)$. Then a continuity argument implies that $h_n\to h$, where $h:[-1,1]\to [-1,1]$ is
a monotone increasing continuous map preserving combinatorics
$$
h:\Bbb{B}_n(g)\to \Bbb{B}_n(\hat{f}).
$$
The partitions $\Bbb{B}_n(\hat{f})$ will get finer and finer. This is a consequence of the fact that $\hat{f}$ does not have homtervals. The map $h$
is continuous. It can happen that $g$ has homtervals. In this case the partitions $\Bbb{B}_n(g)$ will not get finer and finer. As consequence the map $h$ will have intervals which are mapped to points, $h$ is a semiconjugacy from
$g$ to $\hat{f}$. We showed that $f$ and $g$ are essentially conjugated.
\hfill\hfill\qed $\,\,$ (Theorem 1.9)
\enddemo

\bigskip
\centerline{\bf References}
\bigskip

\parindent=40pt
\item{[ACT]} A. Arneodo, P. Coullet, C. Tresser, {\it A possible new mechanism for the onset of turbulence}, Physics Letters, {\bf 81A} n. 4 (1981),197-201.

\item{[D]} A.Dold, {\it Lectures on Algebraic Topology}, Springer Verlag 1972.

\item {[GW]} J. Guckenheimer, R. F. Williams, {\it Structure Stability of 
Lorenz attractors}, Publ. Math. IHES (1979), 59-72.

\item{[HS]} J. H. Hubbard, C. T. Sparrow, {\it The Classification of
 Topologically Expansive Lorenz Maps}, Comm. on Pure and App. Math., {\bf XLIII}, (1990), 431-443.

\item {[L]} E. N. Lorenz, {\it Deterministic non-periodic flow.} J. Atmos. Sci. {\bf 20} (1963), 130-141. 

\item{[M]} R.Ma\~n\'e, {\it Hyperbolicity, Sinks and Measure in One-dimensional Dynamics} Commun. Math. Phys. {\bf 100} (1985), 495-524 and Erratum Commun. Math. Phys. {\bf 112} (1987), 721-724.

\item{[MMMS]} M.Martens, W.de Melo, P.Mendes, S.van Strien {\it On Cherry Flows}, Erg.Th. \& Dyn.Sys. {\bf 10} (1990), 531-554.

\item{[MMS]} M. Martens, W. de Melo, S. van Strien, {\it Julia-Fatou-Sullivan
             Theory for real one-dimensional dynamics}, Acta Math. {\bf 168 }
             (1992) 273-318.

\item{[MS]} W.de Melo, S.van Strien, {\it One-dimensional Dynamics}, Springer Verlag.

\item{[MT]} J. Milnor, W. Thurston, {\it On iterated maps of the interval},
            Springer Lecture Notes in Mathematics {\bf 1342} (1988) 465-563.

\item{[MP]} W.de Melo, C.Pugh, {\it On the $C^1$ Brunovski Hypothesis},
            J.Diff. Equations {\bf 113} (1994) 300-337.

\item {[P]} W. Parry, {\it Symbolic dynamics and transformations of the unit interval.} Trans. Amer. Math. Soc. {\bf 122} (1966), 368-378.

\item{[R]} A.Rovella, {\it The dynamics of perturbations of contracting Lorenz Maps}, Bul.Soc.Brazil. Mat. (N.S.) {\bf 24} (1993) no.2, 233-259.

\bye